\RequirePackage{ifpdf}
\ifpdf 
\documentclass[pdftex]{sigma}
\else
\documentclass{sigma}
\fi

\font\sg=cmbsy10 scaled \magstep1
\newcommand{\scal}{\kern .7mm \hbox{\sg \char'001} \kern .7mm}

\numberwithin{equation}{section} \numberwithin{theorem}{section}
\numberwithin{definition}{section} \numberwithin{remark}{section}
\numberwithin{note}{section}

\begin{document}

\allowdisplaybreaks

\renewcommand{\PaperNumber}{036}

\FirstPageHeading

\renewcommand{\thefootnote}{$\star$}

\ShortArticleName{Dynamical Equations of Non-Holonomic Systems}

\ArticleName{A `User-Friendly' Approach to the Dynamical\\
Equations of  Non-Holonomic Systems\footnote{This paper is a
contribution to the Vadim Kuznetsov Memorial Issue ``Integrable
Systems and Related Topics''. The full collection is available at
\href{http://www.emis.de/journals/SIGMA/kuznetsov.html}{http://www.emis.de/journals/SIGMA/kuznetsov.html}}}

\Author{Sergio BENENTI}

\AuthorNameForHeading{S. Benenti}

\Address{Department of Mathematics, University of Turin, Italy}
\Email{\href{mailto:sergio.benenti@unito.it}{sergio.benenti@unito.it}}
\URLaddress{\url{http://http://www2.dm.unito.it/~benenti/}}

\ArticleDates{Received November 29, 2006, in f\/inal form February
13, 2007; Published online March 01, 2007}

\Abstract{Two ef\/fective methods for writing the dynamical
equations for non-holonomic systems are illustrated. They are
based on the two types of representation of the constraints: by
parametric equations or by implicit equations. They can be applied
to linear as well as to non-linear constraints. Only the basic
notions of vector calculus on Euclidean 3-space and on tangent
bundles are needed. Elementary examples are illustrated.}


\Classification{37J60; 70F25}

\renewcommand{\thefootnote}{\arabic{footnote}}
\setcounter{footnote}{0}

\section{Preamble}

The classical theory of non-holonomic dynamical systems, even in
recent times, is treated in a~growing number of papers. Most of
them require the use of modern dif\/ferentiable and algebraic
structures which are not familiar to non-mathematicians working on
concrete applications. On the other hand, several papers are
dedicated to the analysis of special non-holonomic mechanical
systems, quite interesting but treated with \emph{ad hoc} methods.
These are the reasons why I think useful to propose a \emph{ready
to use} approach to the dynamics of non-holonomic systems,
requiring the knowledge of the basic notions of the vector
calculus on the Euclidean three-space and on tangent bundles only,
and avoiding the use of cotangent bundles (Hamiltonian framework)
and jet-bundles.

In the present paper I illustrate two dif\/ferent, general,
ef\/fective and concise methods for writing the dynamical
equations of a given non-holonomic system. These methods
correspond to the two ways of representing kinematical
constraints: by parametric equations or by implicit equations.
They can be applied to linear as well as to non-linear constraints
and lead to dif\/ferent (but obviously equivalent) dynamical
systems.

In order to make this paper self-contained, a straightforward
approach to the Gauss principle and to the Gibbs--Appell equations
is illustrated. Our starting point will be the well understandable
Newton dynamical equations for a system of massive points.

The tutorial character of this paper does not exclude the presence
of some novelties.

Many articles and books have been consulted in writing this paper.
To cite all them would make the list of references quite long.  On
the other hand, such a long list would not f\/it with the purposes
of this paper. Anyway, I must mention the excellent recent books
\cite{Bullo-Lewis, Monforte},  from which I have got many useful
hints. They are very well readable and advisable to non-experts
which would like to go into this matter in more depth, and to know
its present developments and applications.


\section{Introduction}

Let us consider a mechanical system with a well def\/ined
conf\/iguration manifold $Q_n$. The dimension $n$ of $Q$ is the
number of degrees of freedom of the system. Let us  denote by
$q=(q^i)$ any (Lagrangian) coordinate system on $Q$ and by
$(q,\dot q)=(q^i,\dot q^i)$ the associated coordinates on the
tangent bundle $TQ$. This tangent bundle is the {\bf space of the
kinematical states}.

A {\bf kinematical constraint} is given by a subset $C\subset TQ$.
A special case is that of \emph{regular constraint}:

\begin{definition}
A kinematical constraint $C\subset TQ$ is said to be {\bf regular}
if
\begin{itemize}\itemsep=0pt
\item $C$ is a submanifold of dimension $n+m$, $m<n$; \item for
all $q\in Q$, $F_q=C\cap T_qQ$ is a submanifold of dimension $m$;
\item the restriction to $C$ of the tangent f\/ibration
$\tau_Q\colon TQ\to Q$ is a \emph {surjective
submersion}\footnote{If this notion is not understood by the
reader, he can look at the equivalent conditions (\ref{RC1}) and
(\ref{RC2}) below. The def\/inition of regularity is taken from
\cite{Oliva-Kobayashi} and \cite{Marle}. In \cite{Marle} an
extension of this def\/inition is given: it requires the existence
of a submanifold $Q_1\subset Q$ such that $C$ is a submanifold of
$TQ_1$ and the restriction to $C$ of the tangent f\/ibration is a
submersion. For our purposes we do not need to consider this more
general case.}.
\end{itemize}
The constraint is {\bf linear} if each $F_q$ is a subspace (see
Fig.~1).
\end{definition}

\begin{figure}[ht]
\begin{center}
\includegraphics[width=7cm]{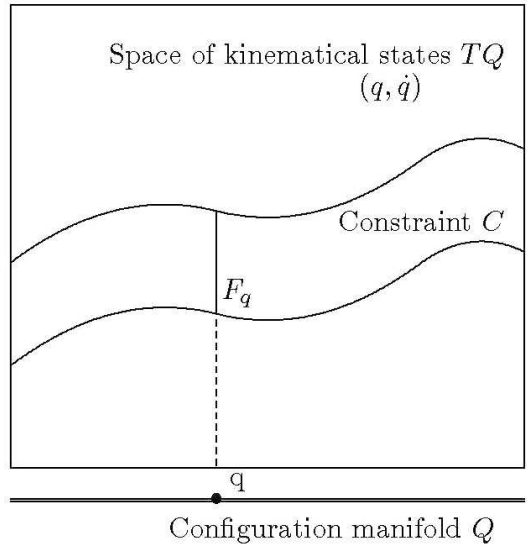}
\caption{A regular constraint is a f\/ibration onto $Q$.}
\end{center}
\end{figure}

A regular constraint can be represented in two ways:
\begin{itemize}\itemsep=0pt
\item {\bf Parametric representation}: it is described by $m$
equations,
  \begin{equation}\label{AA}
  \dot q^i=\psi^i(q,z),
    \end{equation}
where $z=(z^\alpha)$, $\alpha=1,\ldots,m<n$ are called
\emph{parameters}. Note that $(q,z)$ can be interpreted as local
coordinates on $C$. \item {\bf Implicit representation}: it is
described by $r=n-m$ independent equations\footnote{This means
that the dif\/ferentials $dC^a$ are linearly independent at each
point of $C$.},
  \begin{equation}\label{IR}
  C^a(q,\dot q)=0, \qquad a=1,\ldots, r.
    \end{equation}
\end{itemize}

In these two representations, the regularity of the constraint is
represented by conditions
   \begin{equation}\label{RC1}
\hbox{\rm rank}\, [\psi^i_\alpha]_{n\times m}=m, \qquad \hbox{\rm
rank}\, [C^a_i]_{n\times r}=r,
   \end{equation}
respectively, where
   \begin{equation}\label{RC2}
  \psi^i_\alpha\doteq\dfrac{\partial\psi^i}{\partial z^\alpha},
\qquad
  C^a_i\doteq\dfrac{\partial C^a}{\partial \dot q^i}.
   \end{equation}

Note that the regularity conditions may be not fulf\/illed at
certain states, that we call {\bf singular states}. They always
occurs, for instance, for non-linear homogeneous constraints (see
Remark~\ref{r:SS}).

In the following, with the exception of special remarks, the
constraints are assumed to be non-linear. However, in the case of
linear constraints, we shall assume~-- without loss of
generality~-- that the functions $\psi^i(q,z)$ are linear in the
parameters, so that $\psi^i=\psi^i_\alpha(q)\,z^\alpha$, or that
the functions $C^a(q,\dot q)$ are linear in the $\dot q$, so that
$C^a=C^a_i(q)\,\dot q^i$.

The leading idea of the {\bf f\/irst method} is to consider the
parametric equations (\ref{AA})  of the constraint as a f\/irst
set of dynamical equations, to be completed by a second set of
equations of the kind $\dot z^\alpha=Z^\alpha(q,z)$:
  \begin{gather}
\dot q^i=\psi^i(q,z), \nonumber\\
\dot z^\alpha=Z^\alpha(q,z).\label{AB}
  \end{gather}
This is a f\/irst-order system of ODE's. Once the initial
conditions $(q_0,z_0)$ are f\/ixed, they give a~unique actual
motion of the mechanical system. By {\bf actual motion} we mean a
physical motion obeying to the constraints.

The explicit expressions of the functions $Z^\alpha(q,z)$ depend
of course on the given func\-tions $\psi^i(q,z)$.  But they depend
also on the dynamical principles we assume.

The Newton dynamical equation for a material point will be the
only physical principle on which we shall base our approach.
Furthermore, we shall assume that the constraints are {\bf ideal}
or {\bf perfect}, according to a suitable mathematical
def\/inition: this means to accept a certain \emph{constitutive
condition} on the constraint as a postulate. In this way, the
\emph{Gauss} `\emph{principle}' will follow as  a `\emph{theorem}'
from the Newton equations of a system of particles.

The {\bf second method} presented here is based on the implicit
representation of the constraints. It is a revisitation, with
improvements and simplif\/ications, of the well-known
Lagrange-multipliers method. In this context we get the explicit
expression of the reactive forces. This is useful, in the concrete
applications, for measuring the stress that the constraints have
to support.

These two methods lead to dif\/ferent f\/irst-order dynamical
equations i.e., to dif\/ferent vector f\/ields:
\begin{itemize}\itemsep=0pt
\item The f\/irst method (parametric representation)  produces a
vector f\/ield $\mathbf Z$ on the constraint manifold $C$. The
integral curves of $\mathbf Z$ give all actual motions -- see
Fig.~2. \item The second method (implicit representation) produces
a vector f\/ield $\mathbf D$ on the whole $TQ$ but tangent to $C$.
Only its restriction to $C$ has in fact a physical meaning i.e.,
only its integral curves which start from a point of $C$ (and
which will lie on $C$) represent actual motions -- see Fig.~3. As
a byproduct, this method gives the explicit expressions of the
reactive forces, which can be estimated along any actual motion.
\end{itemize}

\begin{figure}[h]
\begin{minipage}{75mm}
\centerline{\includegraphics[width=7cm]{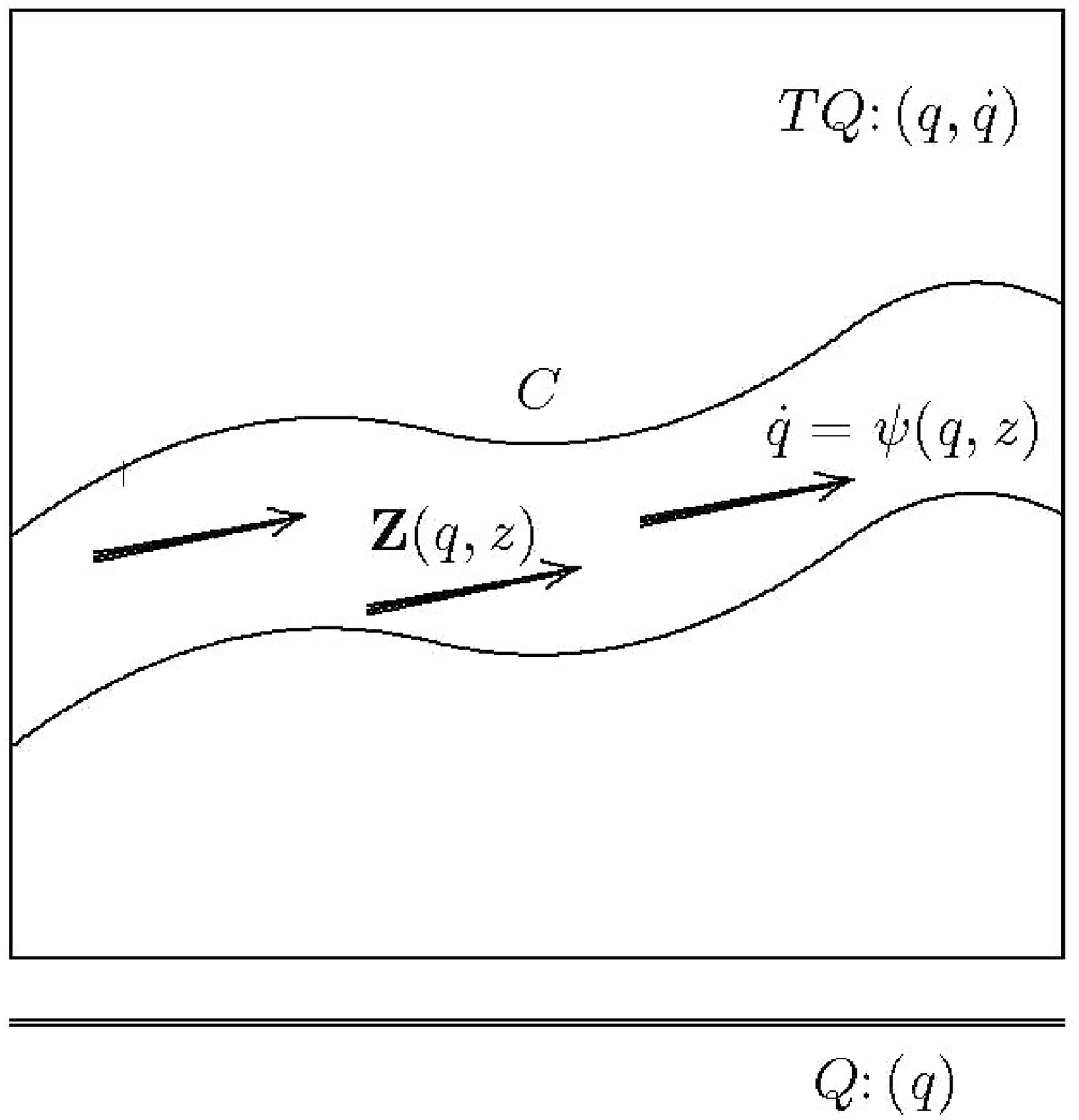}}
\vspace{-2mm} \caption{Parametric representation of $C$ -- The
vector f\/ield $\mathbf Z$.}
\end{minipage}\hfil\quad\hfil
\begin{minipage}{75mm}
\centerline{\includegraphics[width=7cm]{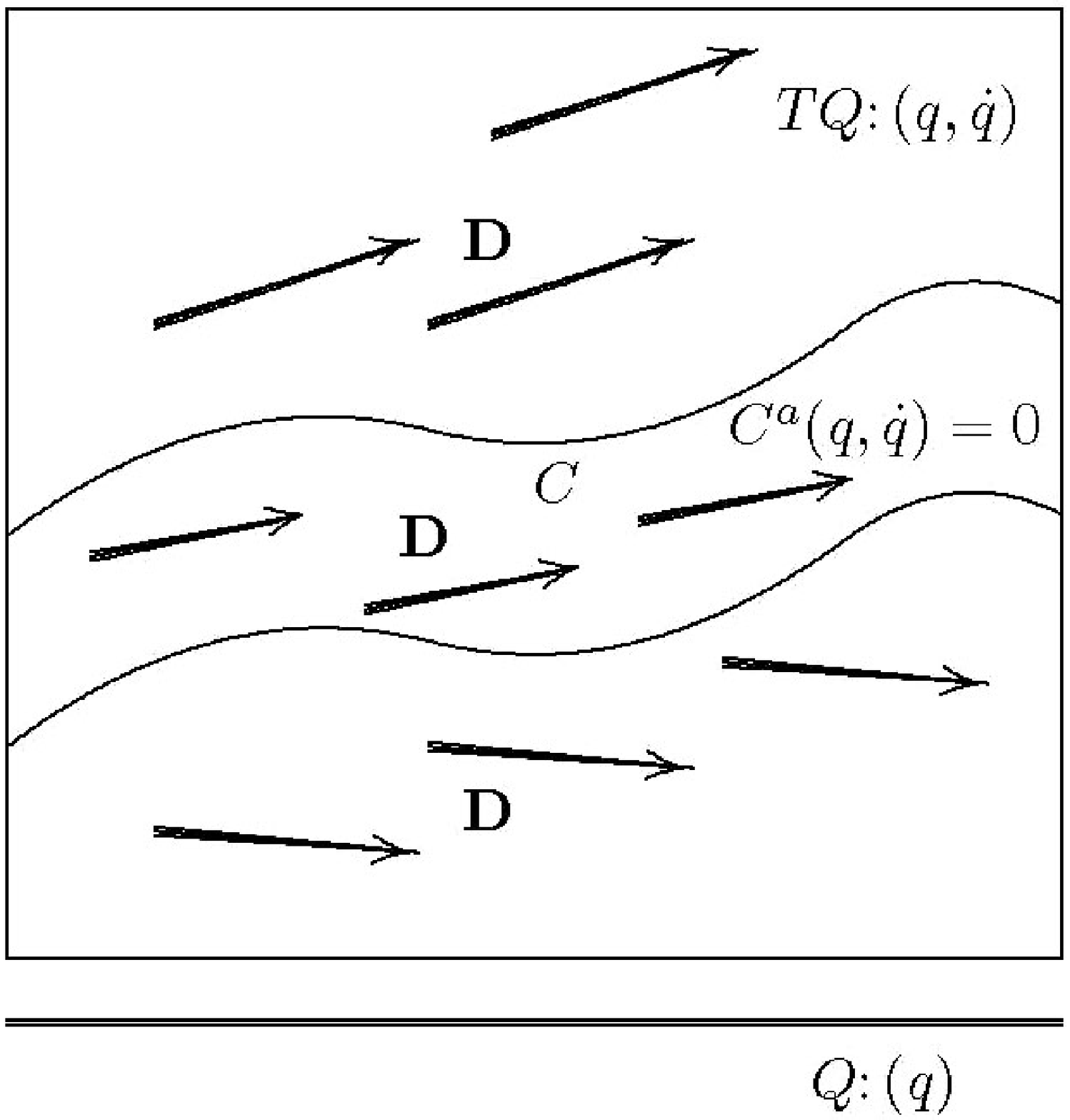}}
\vspace{-2mm} \caption{Implicit representation of $C$ -- The
vector f\/ield $\mathbf D$.}
\end{minipage}
\end{figure}


\section{Ideal constraints}

At a {\it microscopical} level, a mechanical system is made of a
collection of material points $(P_\nu,m_\nu)$, $\nu=1,\ldots,N$.
The position vector $\mathbf r_\nu$ of each point $P_\nu$ in the
Euclidean three-space is a function~$\mathbf r_\nu(q)$ of the
chosen Lagrangian coordinates. A {\bf motion} of the system is
then represented by a~time-parametrized curve $q^i=q^i(t)$ on $Q$,
and at each conf\/iguration $q\in Q$ all possible velocities of
the points are given by
  \begin{equation} \label{VNU}
  \mathbf v_\nu=\dfrac{\partial\mathbf r_\nu}{\partial q^i}\;\dot q^i, \qquad (v^i)\in\mathbb R.
  \end{equation}
A kinematical {\bf state} of the system is the collection of all
pairs position-velocity $(\mathbf r_\nu, \mathbf v_\nu)$ of the
points. The collection of all possible states is then the tangent
bundle $TQ$ of $Q$. At any f\/ixed state the accelerations
$\mathbf a_\nu$ are given by
  \begin{equation}\label{ANU}
\mathbf a_\nu =\dfrac{\partial^2\mathbf r_\nu}{\partial
q^i\partial q^j}\;\dot q^i\dot q^j  +\dfrac{\partial\mathbf
r_\nu}{\partial q^i}\; \frac{d\dot q^i}{dt}.
   \end{equation}
Let us consider the parametric representation (\ref{AA}) of the
constraint. Then the velocity and the acceleration of each point,
compatible with the kinematical constraint, are
\begin{gather*}
\mathbf v_\nu (q,z)=\dfrac{\partial \mathbf r_\nu}{\partial q^i}\,\psi^i,\\
\mathbf a_\nu (q,z,\dot z)=\dfrac{\partial^2\mathbf
r_\nu}{\partial q^i\partial q^j}\,\psi^i\psi^j
+\dfrac{\partial\mathbf r_\nu}{\partial q^i}
\,\bigg(\dfrac{\partial\psi^i}{\partial q^j}\,\psi^j
+\dfrac{\partial\psi^i}{\partial z^\alpha}\,\dot z^\alpha\bigg).
\end{gather*}
Let us write these equations in the more compact form, according
to the above-given notation:
\begin{gather}
\mathbf v_\nu (q,z) =\partial_i\mathbf r_\nu\,\psi^i,\nonumber\\
\mathbf a_\nu (q,z,\dot z) =\partial_{ij}\mathbf
r_\nu\,\psi^i\psi^j +\partial_i\mathbf
r_\nu\,\big(\psi^i_j\,\psi^j+\psi^i_\alpha\,\dot
z^\alpha\big).\label{VA}
\end{gather}
As we shall see, the additional parameters $\dot z$ play a crucial
role. Note that, at any f\/ixed conf\/iguration $q\in Q$, the
parameters $z$ determine all kinematical states $(q,v)$ compatible
with the constraints. Intuitively, the parameters $\dot z$
determine a {\it small displacement} from a state $(q,\dot q)$ to
another close state $(q,\dot q')$ (with the same conf\/iguration
$q$). For a precise def\/inition, let us decompose the
acceleration (\ref{VA}) into the sum
  \begin{equation}
\mathbf a_\nu=\mathbf a_{0\nu}+\mathbf a_{\alpha\nu}\,\dot
z^\alpha, \label{AD}
 \end{equation}
where
  \begin{equation}
\mathbf a_{0\nu}  \doteq\partial_{ij}\mathbf r_\nu\,\psi^i\psi^j
+\partial_i\mathbf r_\nu\,\psi^i_j\,\psi^j, \qquad \mathbf
a_{\alpha\nu}  \doteq \partial_i\mathbf r_\nu\,\psi^i_\alpha.
\label{ADD}
  \end{equation}
We observe that the vectors $\mathbf a_{0\nu}$ depend only on the
state of the system. We say that

\begin{definition}
The  second vector in $(\ref{AD})$
  \begin{equation}
\mathbf w_\nu\doteq\mathbf a_{\alpha\nu}\,\dot z^\alpha,
\label{VD}
  \end{equation}
represents the {\bf virtual displacements} of the point $P_\nu$ at
the given kinematical state.
\end{definition}

\begin{remark}
Note that the parameters $\dot z^\alpha$ span all possible virtual
displacements at a given state. The mechanical meaning of
\emph{virtual displacement}, which is strictly related to that of
\emph{ideal constraint} (Def\/inition \ref{IC}) would require a
detailed discussion. We can skip it simply by accepting
equation~(\ref{VD})  as a \emph{mathematical definition}, since it
will be justif\/ied f\/irst, by the fact that with such a
def\/inition the Gauss principle becomes a consequence of the
Newton dynamical equations (Theorem~\ref{t:GP}) and second, by the
fact that the reactive forces of ideal constraints are not
dissipative (Remark~\ref{r:WR}).

It is customary to associate the intuitive idea of `virtual
displacement' with that of `virtual velocity', as a limit of a
`small' displacement between two conf\/igurations' of the system.
Instead, within the present context, a `virtual displacement' is a
`small' displacement between kinema\-ti\-cal states (conf\/igurations
plus velocities), so it is associated with the intuitive idea of
`virtual acceleration'. This viewpoint is in fact coherent with
the philosophy of the Gauss principle, which deals with
accelerations.
\end{remark}

We assume for the dynamics of each point $(P_\nu,m_\nu)$ the
Newton equation
 \begin{equation}
m_\nu\,\mathbf a_\nu=\mathbf A _\nu+\mathbf R_\nu, \label{Newton}
 \end{equation}
where $\mathbf A_\nu$ is the {\bf active force} (due to external
f\/ields and internal interactions) and  $\mathbf R_\nu$ is the
{\bf reactive force}: it has the role of making the constraint
fulf\/illed.

\begin{remark}
The idea of `reactive force' arises from the Newtonian philosophy,
according which any action deviating a point from the uniform
rectilinear motion (in an inertial reference frame) is a `force',
mathematically represented by a vector. Thus, the presence of a
kinematical constraint must be represented by a vector, called
`reactive force', to be summed to the `active force', which in
turns represents the action of f\/ields present in the space and
independent from the constraints (gravitational,
electromagnetical, centrifugal, Coriolis, etc.).
\end{remark}

\begin{definition}\label{IC}
Non-holonomic constraints are said to be {\bf ideal} or {\bf
perfect} if
  \begin{equation}
\sum_\nu\,\mathbf R_\nu\scal\mathbf w_\nu=0
\label{VW}
  \end{equation}
for all virtual displacements $\mathbf w_\nu$.
\end{definition}

We consider equation~(\ref{VW}) as a {\bf constitutive condition}
for the constraint: it says which kind of reactive forces the
constraint is able to supply in order to be satisf\/ied along any
motion. It is straightforward to check that for linear
constraints, or simply for holonomic constraints (which do not
involve velocities), equation (\ref{VW}) reduces to the classical
{\bf virtual work principle}.
 The validity of such a constitutive condition is a matter of  theoretical and experimental analysis of the behavior of the mechanical system one is dealing with.

\begin{theorem}\label{t:IC}
Let
\begin{equation}\label{RI}
R_i\doteq \sum_\nu\;\mathbf R_\nu\scal\dfrac{\partial\mathbf
r_\nu}{\partial q^i},
  \end{equation}
be the {\bf Lagrangian reactive forces}. Then the definition {\rm
(\ref{VW})} of ideal constraint is equivalent to the following
equations,
  \begin{gather}
R_i\;\psi^i_\alpha=0, \label{RP}
\\
\label{RIM} R_i=\lambda_a\,C^a_i,
\end{gather}
in the parametric and in the implicit representation,
respectively.
\end{theorem}

Equation (\ref{RIM}) means that the components $R_i$ are linear
combinations of the functions $C^a_i(q,\dot q)$.

\begin{proof}
Put equation (\ref{VD}) and (\ref{ADD}) in equation  (\ref{VW}),
 \[
0= \sum_\nu \mathbf R_\nu\scal\mathbf a_{\alpha\nu}\,\dot
z^\alpha=
 \sum_\nu \mathbf R_\nu \scal \partial_i\mathbf r_\nu\,\psi^i_\alpha\,\dot z^\alpha=R_i\,\psi^i_\alpha\,\dot z^\alpha.
 \]
This proves equation (\ref{RP}). By dif\/ferentiating the identity
$C^a(q,\psi(q,z))=0$, we get the following relations between the
two representations,
  \begin{equation}
C^a_i\,\psi^i_\alpha=0, \qquad
\partial_iC^a+C^a_j\,\psi^j_i=0,
\label{ID}
  \end{equation}
  where
  \begin{equation}
C^a_i\doteq \dfrac{\partial C^a}{\partial\dot q^i}, \qquad
\psi^j_i\doteq \partial_i\psi^j,\qquad
\partial_i\doteq\dfrac{\partial}{\partial q^i}.
 \label{ID1}
\end{equation}
Due to the f\/irst identity (\ref{ID}), equation (\ref{RP}) is
then equivalent to (\ref{RIM})
\end{proof}

\begin{remark} \label{r:WR}
Equations $(\ref{ID1})_1$ and (\ref{RIM}) show that the ideal
constraints do not dissipate ener\-gy: the power of the reactive
forces is zero (we are dealing with time-independent constraints).
See \cite{Caratheodory}.
\end{remark}

\begin{remark} \label{r:LE}
By a well-known process, we pass from the microscopical level to
the `\emph{macrosco\-pical}' one i.e., to the Lagrange equations
\begin{equation}\label{LE}
\dfrac{d}{dt}\left(\dfrac{\partial K}{\partial \dot q^i}\right)-
\dfrac{\partial K}{\partial q^i}=A_i+R_i,
\end{equation}
where
  \begin{equation}\label{K}
  K=\frac 12\;g_{ij}\,\dot q^i\dot q^j, \qquad g_{ij}\doteq  \sum_\nu
  m_\nu\;\partial_i\mathbf r_\nu\scal\partial_j\mathbf r_\nu.
   \end{equation}
At the right hand side of the Lagrange equations we f\/ind the sum
of the {\bf active Lagrangian forces}
\begin{equation}\label{AI}
A_i\doteq \sum_\nu\;\mathbf A_\nu\scal\dfrac{\partial\mathbf
r_\nu}{\partial q^i},
  \end{equation}
and the reactive Lagrangian forces (\ref{RI}).
\end{remark}


\section{The Gauss principle}\label{GP}

At the microscopical level we introduce the quantity
  \begin{equation}
G\doteq {\frac 12}\, \sum_\nu m_\nu\bigg(\mathbf
a_\nu-\dfrac{\mathbf A_\nu}{m_\nu}\bigg)^2. \label {G1}
 \end{equation}
The active forces $\mathbf A_\nu$ are known functions of the state
$(q,\dot q)$. Thus, due to the parametric equations of the
constraints, $G$ becomes a function of $(q,z)$. Moreover, even for
active forces depending on the velocities, $\mathbf A_\nu$ does
not depend on $\dot z$,
  \[
\dfrac{\partial \mathbf A_\nu}{\partial \dot z^\alpha}=0.
  \]
Thus, due to equations (\ref{AD}) and (\ref{ADD}), along any
motion satisfying the constraints  we have
  \begin{equation}\label{GZ}
 \dfrac{\partial G}{\partial \dot z^\alpha}
= \sum_\nu m_\nu \left(\mathbf a_\nu-\dfrac{\mathbf
A_\nu}{m_\nu}\right)\scal \dfrac{\partial \mathbf a_\nu}{\partial
\dot z^\alpha} = \sum_\nu m_\nu \left(\mathbf a_\nu-\dfrac{\mathbf
A_\nu}{m_\nu}\right)\scal \mathbf a_{\alpha\nu}.
  \end{equation}
Furthermore,
  \[
\dfrac{\partial^2 G}{\partial \dot z^\alpha \partial \dot z^\beta}
=  \sum_\nu m_\nu\;\dfrac{\partial \mathbf a_\nu}{\partial \dot
z^\beta}\scal\mathbf a_{\alpha\nu} = \sum_\nu m_\nu\;\mathbf
a_{\beta\nu}\scal \mathbf a_{\alpha\nu}  = \sum_\nu m_\nu\;
\psi^i_\alpha\,\psi^j_\beta\;\partial_i\mathbf
r_\nu\scal\partial_j\mathbf
r_\nu=g_{ij}\,\psi^i_\alpha\,\psi^j_\beta.
  \]
Then, if we introduce the functions
  \begin{equation}\label{Gcov}
G_{\alpha\beta}\doteq  g_{ij}\,\psi^i_\alpha\,\psi^j_\beta,
  \end{equation}
we get
  \begin{equation}\label{GG}
\dfrac{\partial^2 G}{\partial \dot z^\alpha \partial \dot
z^\beta}=G_{\alpha\beta}.
  \end{equation}
Since the matrix $[\psi^i_\alpha]$ has maximal rank, the symmetric
matrix $[G_{\alpha\beta}]$ is regular and positive-def\/inite as
well as $[g_{ij}]$.

\begin{theorem}\label{t:GP}
Assume the Newton equations $m_\nu\,\mathbf a_\nu=\mathbf
A_\nu+\mathbf R_\nu$ for each point $P_\nu$. Then, at any state of
any actual motion the quantity $G$ takes a minimal value  \rm
(Gauss principle) \it if and only if the constraints are ideal.
\end{theorem}

\begin{proof}
Write the Newton equations in the form
  \[
m_\nu\,\left(\mathbf a_\nu-\dfrac{\mathbf
A_\nu}{m_\nu}\right)=\mathbf R_\nu.
  \]
Then, due to equations (\ref{AD}), (\ref{ADD}) and (\ref{RI}),
\begin{equation}\label{GZR}
\dfrac{\partial G}{\partial \dot z^\alpha}= \sum_\nu \mathbf
R_\nu\scal\dfrac{\partial\mathbf a_\nu}{\partial \dot z^\alpha}=
\sum_\nu \mathbf R_\nu\scal\mathbf a_{\alpha\nu}
=R_i\,\psi^i_\alpha.
  \end{equation}
This shows that for ideal constraints, see equation (\ref{RP}),
the Newton equations imply
  \begin{equation}\label{GZ0}
\dfrac{\partial G}{\partial \dot z^\alpha}=0,
  \end{equation}
at any state along any actual motion. Due to equation (\ref{GG}),
being $[G_{\alpha\beta}]$ positive, at the statio\-nary states for
which equation (\ref{GZ0}) holds, the function $G$ has a strong
minimum. (ii)~Conversely, assume that the Gauss principle holds
true. Then equation (\ref{GZ0}) is satisf\/ied, so that from
(\ref{GZR}) we get $R_i\,\psi^i_\alpha=0$. This means that the
constraint is ideal (Theorem~\ref{t:IC}).
\end{proof}

\begin{remark}
The vector $\mathbf A_\nu/m_\nu$ is the acceleration of the point
$P_\nu$ in a {\bf free motion}, free from the constraints. Let us
denote it by $\mathbf a_\nu^f$. As a consequence, the function $G$
can be also def\/ined as
  \begin{equation}
G\doteq {\frac 12}\,\sum_\nu m_\nu\big(\mathbf a_\nu-\mathbf
a_\nu^f\big)^2, \label{G2}
  \end{equation}
and Theorem \ref{t:GP} can be reformulated as follows:
\end{remark}

\begin{theorem}\label{TH:2}
Let $\mathbf r_\nu(t)$ and $\mathbf r_\nu^f(t)$  be \underbar{two}
motions of the system $P_\nu$ such that for $t=t_0$ the
corresponding states coincide i.e.,
  \[
\mathbf r_\nu(t_0)=\mathbf r_\nu^f(t_0), \qquad
 \mathbf v_\nu(t_0)=\mathbf v_\nu^f(t_0).
  \]
Assume that $\mathbf r_\nu^f(t)$ is a free motion. Then, at this
state, and for any motion compatible with ideal constraints, the
actual accelerations $\mathbf a_\nu(t_0)$, are such that $G$ takes
a minimal value.
\end{theorem}

\begin{remark}
For any arbitrary motion,
  \[
\mathbf v_\nu=\dfrac{\partial \mathbf r_\nu}{\partial q^i}\,\dot
q^i,\qquad \mathbf a_\nu=\dfrac{\partial\mathbf r_\nu}{\partial
q^i\partial q^j}\,\dot q^i\dot q^j +\dfrac{\partial\mathbf
r_\nu}{\partial q^i}\,\ddot q^i.
  \]
Then, at any \underbar{f\/ixed} \underbar{state} compatible with
the constraints we have
  \[
\mathbf a_\nu-\mathbf a_\nu^f=\dfrac{\partial\mathbf
r_\nu}{\partial q^i}\,(\ddot q^i-\ddot q^i_f).
  \]
Due to the def\/inition (\ref{K}) of $g_{ij}$, from (\ref{G2}) we
get the so-called {\bf Lipschitz expression} of the function $G$:
  \begin{equation}
G={\frac 12}\, g_{ij}\,(\ddot q^i-\ddot q^i_f)(\ddot q^j-\ddot
q^j_f).
  \label{G3}
 \end{equation}
Note that in this expression the Christof\/fel symbols are not
involved.
\end{remark}

\section[The Gibbs-Appell equations]{The Gibbs--Appell equations}

Let us go back to the def\/inition (\ref{G1}) of the function $G$.
If we introduce the functions
  \begin{equation}
S\doteq {\frac 12} \sum_\nu m_\nu\,\mathbf a_\nu^2
,\qquad S_1\doteq {\frac 12}\sum_\nu\frac
1{m_\nu}\,\mathbf A_\nu^2, \qquad S_2\doteq
\sum_\nu\mathbf A_\nu\scal\mathbf a_\nu, \label{ABC}
\end{equation}
then we have the decomposition
  \[
G=S+S_1-S_2.
  \]
The function $S$ is called the {\bf energy of the accelerations}.
We observe that
  \[
\frac{\partial S_1}{\partial \dot z^\alpha}=0
  \]
and that, due to the second equation (\ref{VA}) and the
def\/inition (\ref{AI}) of active Lagrangian force,
  \[
\frac{\partial S_2}{\partial \dot z^\alpha}=\sum_\nu\mathbf
A_\nu\scal\partial_i\mathbf
r_\nu\,\psi^i_\alpha=A_i\,\psi^i_\alpha.
  \]
Thus,
  \[
\frac{\partial G}{\partial \dot z^\alpha}=\frac{\partial
S}{\partial \dot z^\alpha}-A_i\,\psi^i_\alpha.
  \]
Due to the Gauss principle (Theorem \ref{t:GP}), this proves

\begin{theorem}
The Gauss principle is equivalent to equations
 \begin{equation}
\frac{\partial S}{\partial \dot z^\alpha}=A_\alpha, \label{GT1}
  \end{equation}
where the function
  \begin{equation}
 S(q,z,\dot z)\doteq {\frac 12} \sum_\nu m_\nu\,\mathbf a_\nu^2
 \label{GT2}
  \end{equation}
is determined by the expression \rm (\ref{VA}) \it of the
accelerations, and
 \begin{equation}
A_\alpha\doteq A_i\psi^i_\alpha. \label{GT3}
 \end{equation}
\end{theorem}

\begin{remark}
Equations (\ref{GT1}) are the celebrated {\bf Gibbs--Appell
equations}. The equivalence between these equations and the Gauss
principle is highlighted within the framework of the parametric
representation (\ref{AA}) of the constraints.
\end{remark}

\begin{remark}
The quantities $A_\alpha$ can be computed by writing the the {\bf
virtual power} of the active forces:
  \begin{equation}
W \doteq \sum_\nu\,\mathbf A_\nu\scal\mathbf w_\nu =\sum_\nu
\mathbf A_\nu\scal\partial_i\mathbf r_\nu\,\psi^i_\alpha\,\dot
z^\alpha= A_i \,\psi^i_\alpha\,\dot z^\alpha= A_\alpha\,\dot
z^\alpha.
  \end{equation}
\end{remark}

\section[The explicit form of the Gibbs-Appell equations]{The explicit form of the Gibbs--Appell equations}

Both sides of the Gibbs--Appell equations (\ref{GT1}) are
functions of $(q,z,\dot z)$. Let us solve them w.r.to the
variables $\dot z^\alpha$. To this end, it is crucial to observe
that by using the second equations (\ref{VA}) we get for the
function $S$ (\ref{GT2}) the expression
  \[
S={\frac 12}\, g_{ij}\,\psi^i_\alpha\psi^j_\beta\,\dot
z^\alpha\dot z^\beta+ \sum_\nu m_\nu\,\partial_{ij}\mathbf
r_\nu\scal\partial_k \mathbf
r_\nu\;\psi^i\psi^j\psi^k_\alpha\,\dot z^\alpha+S_0,
 \]
where $S_0$ is a function dependent on $(q,z)$ only. Then, this
function is not involved by the Gibbs--Appell equations and $S$
can be replaced by
 \begin{equation}
S_*={\frac 12}\, g_{ij}\,\psi^i_\alpha\psi^j_\beta\,\dot
z^\alpha\dot z^\beta+ \sum_\nu m_\nu\,\partial_{ij}\mathbf
r_\nu\scal\partial_k \mathbf
r_\nu\;\psi^i\psi^j\psi^k_\alpha\,\dot z^\alpha. \label{SS}
 \end{equation}
Now we show that this new function $S_*$ assumes a very
interesting expression. Let us introduce the functions
 \[
\xi_{ijk}(q)\doteq \sum_\nu m_\nu\,\partial_{ij}\mathbf
r_\nu\scal\partial_k \mathbf r_\nu.
  \]
Since
  \[
  \xi_{ijk}= \sum_\nu m_\nu\,\partial_i(\partial_j\mathbf r_\nu\scal\partial_k\mathbf r_\nu)- \sum_\nu m_\nu\,(\partial_j\mathbf r_\nu\scal\partial_{ik}\mathbf r_\nu)=\partial_ig_{jk}-\xi_{ikj},
  \]
by a cyclic permutation of the indices we get
  \[
\xi_{ijk}+\xi_{ikj}=\partial_ig_{jk},\qquad
\xi_{jki}+\xi_{jik}=\partial_jg_{ki},\qquad
\xi_{kij}+\xi_{kji}=\partial_kg_{ij}.
  \]
By summing the f\/irst two equations and subtracting the third
one, since $\xi_{ijk}$ is symmetric in the f\/irst two indices, we
get $\xi_{ijk}=\Gamma_{ijk}$, where
  \[
\Gamma_{ijk}\doteq{\frac
12}\,(\partial_ig_{jk}+\partial_jg_{ki}-\partial_kg_{ij})
  \]
are the Christof\/fel symbols of the metric tensor $g_{ij}$ (the
coef\/f\/icients of the Levi-Civita connection). As a consequence,
if we recall the def\/inition (\ref{Gcov}) of $G_{\alpha\beta}$,
the function $S_*$ (\ref{SS}) can be written as
 \[
S_*={\frac 12}\, G_{\alpha\beta}\,\dot z^\alpha\dot
z^\beta+\Gamma_{ijk}\;\psi^i\psi^j\psi^k_\alpha\,\dot z^\alpha,
 \]
and the Gibbs--Appell equations (\ref{GT1}) assume the form
  \begin{equation}\label{GE}
G_{\alpha\beta}\,\dot
z^\beta+\Gamma_{ijk}\;\psi^i\psi^j\psi^k_\alpha=A_\alpha.
  \end{equation}
Then we can prove
\begin{theorem}
The Gibbs--Appell equations $(\ref{GT1})$ are equivalent to
equations
\begin{equation}\label{NF}
\dot
z^\alpha=G^{\alpha\beta}(A_\beta-\Gamma_{ijk}\;\psi^i\psi^j\psi^k_\beta),
  \end{equation}
where $[G^{\alpha\beta}]$ the inverse matrix of
$[G_{\alpha\beta}]$.
\end{theorem}
\begin{proof}
Indeed, as we remarked in Section \ref{GP}, the matrix
$[G_{\alpha\beta}]$ is regular. If we apply the inverse matrix
$[G^{\alpha\beta}]$ to equations (\ref{GE}), then we get equations
(\ref{NF}).
\end{proof}
Equations (\ref{NF}) are the {\bf explicit form} (or {\bf normal
form}) of the Gibbs--Appell equations~(\ref{GT1}).

\section[The dynamical equations of the first kind]{The dynamical equations of the f\/irst kind}

By setting
\[
\dot z^\alpha=\frac{dz^\alpha}{dt},
\]
equations (\ref{NF}) together with the constraint equations
(\ref{AA}) build up a  f\/irst-order dif\/ferential system, in
normal form, in the unknown functions $q^i(t)$ and $z^\alpha(t)$:
  \begin{gather}
 \dfrac{dq^i}{dt}=\psi^i(q,z),\nonumber \\
\dfrac{dz^\alpha}{dt}=G^{\alpha\beta}\,
(A_\beta-\Gamma_{ijk}\;\psi^i\psi^j\psi^k_\beta). \label{ZZ1}
  \end{gather}
Hence, we can summarize the results so far obtained in the
following
\begin{theorem}  \label{t:Z}
Let $Q$ be the configuration $n$-manifold of a mechanical system,
with local Lagrangian coordinates $q=(q^i)$. Let $TQ$ be the
tangent bundle of $Q$ with canonical coordinates $(q,\dot
q)=(q^i,\dot q^i)$. Let $C\subset TQ$ be a submanifold representing
kinematical time-independent constraints. Assume that the
constraint submanifold $C$ is locally described by parametric
equations $\dot q^i=\psi^i(q,z)$ with $m<n$ parameters
$z=(z^\alpha)$, such that the $n\times m$ matrix
\[
[\psi^i_\alpha]\doteq\left[\dfrac{\partial \psi^i}{\partial
z^\alpha}     \right]
\]
has maximal rank $m$. If the constraints are ideal, then the
actual motions are represented  by functions $q^i=q^i(t)$
satisfying the differential system $(\ref{ZZ1})$\footnote{Any
solution of the dif\/ferential system (\ref{ZZ1}) is of course a
set of functions $q^i(t)$ and $z^\alpha(t)$. But, after the
integration, we can get rid of the functions $z^\alpha(t)$.},
where $\Gamma_{ijk}$ are the Christoffel symbols of the metric
tensor $g_{ij}$ associated with the kinetic energy $K=\tfrac
12\,g_{ij}\,\dot q^i\dot q^j$ and $A_\alpha=A_i\;\psi^i_\alpha$,
where $A_i$ are the Lagrangian active forces.
\end{theorem}

\begin{remark}
The dynamical system (\ref{ZZ1}) is the f\/irst-order system
associated with the vector f\/ield $\mathbf Z$ on the constraint
submanifold $C$, whose components $(Z^i,Z^\alpha)$ w.r.to the
coordinates $(q,z)$ are given by
\begin{gather}
Z^i=\psi^i(q,z),\nonumber \\
Z^\alpha=G^{\alpha\beta}\,
(A_\beta-\Gamma_{ijk}\;\psi^i\psi^j\psi^k_\beta). \label{ZZ2}
  \end{gather}
The $Z^\alpha$ are the `vertical components' of $\mathbf Z$ w.r.to
the projection onto $Q$. The actual motions are the projections
onto $Q$ of the integral curves of $\mathbf Z$.
\end{remark}

\begin{remark}\label{r:Z}
This theorem provides a f\/irst `\emph{recipe}' for writing the
dynamical equations for non-holonomic systems with linear or
non-linear ideal constraints:
\begin{enumerate}\itemsep=0pt
\item Choose Lagrangian coordinates $(q^i)$, write the kinetic
energy of the system $K=\tfrac 12\; g_{ij}\,\dot q^i\dot q^j$, and
extract the $n\times n$ matrix $[g_{ij}]$. \item Choose parametric
equations $\dot q^i=\psi^i(q,z)$ of the constraints, compute the
$m\times n$ matrix~$[\psi^i_\alpha]$, and check its rank. If the
constraints are initially expressed by implicit equations, then
use (for instance) the method illustrated in Remark \ref{r:p}
below for f\/inding parametric equations. \item Compute the
$m\times m$ matrix
$[G_{\alpha\beta}]=[g_{ij}\psi^i_\alpha\psi^j_\beta]$, and the
inverse matrix
$[G^{\alpha\beta}]=[G_{\alpha\beta}]^{-1}$.\footnote{Note that
according to this recipe we do not have to compute the inverse
$n\times n$-matrix of $[g_{ij}]$, but only the inverse of
$[G_{\alpha\beta}]$, whose dimension is $m<n$.} \item Write the
Lagrange equations of the free motions \emph{(i.e., with only
active forces $A_i$)} in the form
\begin{equation} \label {LI}
g_{ij}\,\ddot q^j=L_i(q,\dot q)
  \end{equation}
(note that the formal expression of $L_i$ is $L_i=A_i(q,\dot
q)-\Gamma_{hki}\,\dot q^h\,\dot q^k$) and compute
 \begin{equation}\label{ZI}
\bar Z_i(q,z)=
L_i(q,\psi)=A_i(q,\psi)-\Gamma_{hki}(q)\,\psi^h\psi^k.
 \end{equation}
\item Compute $Z_\alpha=\bar Z_i\psi^i_\alpha$ and
$Z^\alpha=G^{\alpha\beta}\,Z_\beta$. \item Write the dynamical
system
\begin{gather}
 \dfrac{dq^i}{dt}=\psi^i(q,z), \nonumber\\
\dfrac{dz^\alpha}{dt}=Z^\alpha(q,z). \label{ZZR}
  \end{gather}
\end{enumerate}
 \end{remark}

\begin{remark}\label{r:p}
When the constraint submanifold $C\subset TQ$ is described by a
system of implicit independent equations of the kind
 \begin{equation}\label{Ca}
 C^a(q,\dot q)=0, \qquad a=1,\ldots, r, \qquad r=n-m,
 \end{equation}
then we have to transform these equations into parametric
equations. The choice of the parameters $z^\alpha$ is completely
free and it is only a matter of convenience, depending on the
explicit concrete form of the dynamical equations (\ref{ZZR}) we
get. Anyway, since the matrix $[\psi^i_\alpha]$ does not have  the
maxi\-mal rank, equations (\ref{Ca}) can be solved w.r.to $m$ of
the $n$ Lagrangian velocities $\dot q^i$, say~-- up to a
reordering -- w.r.to $\dot q^\alpha$, $\alpha=1,\ldots, m$. This
process leads to considering as parameters these Lagrangian
velocities: $z^\alpha=\dot q^\alpha$. It works very well for
linear or af\/f\/ine constraints, where equations~(\ref{Ca}) have
the form
  \[
C^a_i(q)\,\dot q^i-b_i(q)=0.
  \]
Another possible choice of the parameters, for the linear
constraints only, is that related to the use of {\it
quasi-velocities} or {\it quasi-coordinates} -- see for instance
\cite{Neimark}.
\end{remark}

\begin{remark}
For the analysis of the qualitative (or quantitative) behavior of
a non-holonomic mechanical system (like stability, equilibrium
states, small  oscillations, numerical integration, etc.) we can
apply to $\mathbf Z$ all the known theorems about dynamical
systems.  For instance, a {\bf f\/irst integral} is a function
$F(q,z)$ such that
\begin{equation}\label{FIS}
\psi^i(q,z)\; \dfrac{\partial F}{\partial
q^i}+Z^\alpha(q,z)\;\dfrac{\partial F}{\partial z^\alpha }=0.
  \end{equation}
\end{remark}

\begin{remark}\label{r:SS}
The singular points of the dynamical system (\ref{ZZ1}) are the
solutions $(q,z)$  of the simultaneous equations
\begin{gather}
\psi^i(q,z)=0,\nonumber\\
Z^\alpha(q,z)\doteq \bar Z_i\,\psi^i_\alpha=0.\label{SP}
  \end{gather}
In the case of a homogeneous quadratic constraint,
$\psi^i(q,z)=\tfrac 12 \,\psi^i_{\alpha\beta}(q)\;z^\alpha
z^\beta$, equations~(\ref{SP}) become
\begin{gather}
\psi^i_{\alpha\beta}(q)\,z^\alpha z^\beta=0,\nonumber\\
\bar Z_i\,\psi^i_{\alpha\beta}\,z^\beta=0,\label{SP2}
\end{gather}
being $\psi^i_\alpha=\psi^i_{\alpha\beta}(q)\,z^\beta$. This shows
that, whatever $q$ and $\bar Z_i$, singular points are given by
$z^\alpha=0$. But for $z^\alpha=0$ the matrix $[\psi^i_\alpha]$
does not have the maximal rank, since all its elements vanish.
Hence, at these singular points the constraint $C$ is not regular.
\end{remark}

\begin{remark}
The geometrical picture of the above results gives an intrinsic
meaning of the objects we have introduced. Any vector $\mathbf
v\in TC$ can be represented by a sum
  \[
\mathbf v=v^i\,\partial_i+v^\alpha\,\partial_\alpha,
  \]
where $\partial_i=\partial/\partial q^i$ and
$\partial_\alpha=\partial/\partial z^\alpha$ are interpreted as
pointwise independent vector f\/ields on~$TC$: at each point $x$
of $C$ they span the tangent space $T_xC$. The vectors
$\partial_\alpha$ are `\emph{vertical}' i.e., they are tangent to
the f\/ibers $F_q$ of $C$. Hence, we call~$v^\alpha$ the {\bf
vertical components} of $\mathbf v$, while we call~$v^i$ the {\bf
basic components}. For instance, the basic components of $\mathbf
Z$ are $Z^i=\psi^i(q,z)$ and the vertical components are
$Z^\alpha$. The functions $\psi^i_\alpha$ have the role of
transforming \emph{Latin components}, labeled by Latin indices
$h,i,j,k,\ldots$, into \emph{Greek components}, labeled by Greek
indices $\alpha,\beta,\ldots$. For instance, when it is applied to
a one-form (covariant vector) $\bar Z_i\,dq^i$, we get a
\emph{vertical one form}~$Z_\alpha\,dz^\alpha$, and when it is
applied to the covariant metric tensor $g_{ij}$, we get a~metric
tensor~$G_{\alpha\beta}$ on the f\/ibers of $C$, so that, by
raising the indices of $Z_\alpha$ by the contravariant
metric~$G^{\alpha\beta}$, we get a~vertical vector f\/ield
$Z^\alpha\,\partial_\alpha$.
\end{remark}

\begin{note}
I did not f\/ind equations (\ref{ZZ1}) in the recent and old
articles I have consulted. In fact, it is rather surprising that
the simple idea of considering the parametric representation of
the constraints as a part of the dynamical equations  does not
appear in the major textbooks and treatises on non-holonomic
mechanics. Only recently this idea appeared in a paper of Massa
and Pagani \cite{Massa-Pagani-1991}. Their  general approach,
which is based on the jet-bundle theory and deals with
time-dependent constraints, leads to the introduction of the
vector f\/ield $\mathbf Z$. The elementary approach presented here
is of course quite dif\/ferent and leads, for instance, to
dif\/ferent expressions of the vertical part of $\mathbf Z$. Our
equations (\ref{ZZ1}) should be compared with equations~(3.5b) and
(3.15) of~\cite{Massa-Pagani-1991}. The dif\/ference is that the
second equation (\ref{ZZ1}) is written in terms of the Euclidean
vectors $\mathbf F_\nu$, $\mathbf v_\nu$, while equation~(3.15) of
\cite{Massa-Pagani-1991} $\mathbf Z$ is written in terms of the
Lagrange equations, but still in an implicit form.
\end{note}

\begin{remark}
For linear constraints, the approach presented here is more
general than that of \v Caplygin -- see \cite{Neimark}, Ch.~III,
\S~3, where the coordinates $q^i$ are divided into two groups, say
$(q^a,q^\alpha)$, with $a=1,\ldots,m$ and $\alpha=m+1,\ldots,n$.
The constraint equations are assumed to be of the form
  \[
\dot q^\alpha=\sum_a b^\alpha_a\;\dot q^a,
  \]
where the coef\/f\/icients $b^\alpha_a$ and the Lagrangian $L$ are
assumed to be independent from the coordinates
$(q^a)$.\footnote{In the history of the non-holonomic systems we
can f\/ind the famous equations of Maggi, Volterra, Voronec and \v
Caplygin, dealing with linear constraints. The comparison of these
equations with our approach is left to the reader, who can f\/ind
a detailed discussion in in  \cite{Neimark}, Ch.~3. A neat
illustration of the non-holonomic dynamical equations, with the
essential classical and recent bibliography, can be found in the
book \cite{Bullo-Lewis}.}
\end{remark}


\section{The dynamical equations of the second kind}\label{LM}

About the method for writing the dynamical equations of a
non-holonomic mechanical system so far illustrated   two remarks
are in order:
\begin{itemize}\itemsep=0pt
\item It lies on a parametric representation of the constraint $C$
(however, the vector f\/ield $\mathbf Z$ does not depend on the
chosen parametrization). \item It does not give any information
about the reactive forces.
\end{itemize}

Here, we propose an alternative method for writing the dynamical
equations which is based on any implicit representation of $C$ by
a system of independent equations
  \[
C^a(q,\dot q)=0,
  \]
and which provides a way for  evaluating the reactive forces.

The Lagrange equations (\ref{LE}) -- see Remark \ref{r:LE} --  are
equivalent to the dynamical system
 \begin{equation}\label{XL}
 \mathbf X_\lambda=
\left\{
 \begin{split}
& \dfrac{dq^i}{dt}=\dot q^i, \\
& \dfrac{d\dot q^i}{dt}={}-\Gamma_{hk}^i\,\dot q^h\,\dot
q^k+A^i+R^i_\lambda,
\end{split}
\right.
  \end{equation}
on the tangent bundle $TQ$ of the conf\/iguration manifold $Q$.
Here, $A^i=g^{ij}A_j$ and $R^i_\lambda=g^{ij}\,\lambda_a\,C^a_j$
are the contravariant components of the Lagrangian active and
reactive forces, respectively. The label $\lambda$ points out that
the reactive Lagrangian forces depend on the \emph{a priori}
unknown Lagrangian multipliers $\lambda=(\lambda_a)$, according to
Theorem \ref{t:IC} and equation~(\ref{RIM}).

For a better understanding of what we are going to do, it is
useful to consider the following

\begin{definition}\label{d:VV}
We say that a vector $\mathbf V$ on $TQ$ is {\bf vertical} if it
is tangent, at each point where it is def\/ined,  to the
corresponding f\/iber of $TQ$. This is equivalent to say that it
has the form
  \begin{equation}\label{VV}
  \mathbf V=V^i\,\dfrac{\partial}{\partial \dot q^i}.
 \end{equation}
\end{definition}

As a consequence, the vertical part of the vector $\mathbf
X_\lambda$ is given by
\[
\left(A^i+R^i_\lambda-\Gamma_{hk}^i\,\dot q^h\,\dot q^k\right)
\dfrac{\partial}{\partial \dot q^i},
\]
where the active and reactive forces are represented by the
vertical vectors
  \[
  \mathbf A=A^i\,\dfrac{\partial}{\partial \dot q^i}, \qquad  \mathbf R_\lambda=R^i_\lambda\,\dfrac{\partial}{\partial \dot q^i},
  \]
Hence, the vector $\mathbf X_\lambda$ is decomposed into the sum
 \begin{equation}\label{totalfield}
  \mathbf X_\lambda\doteq\mathbf X_G+\mathbf A+\mathbf R_\lambda,
    \end{equation}
where $\mathbf X_G$ represents the {\bf geodesic f\/low},
\begin{equation}
 \mathbf X_G=
\left\{
 \begin{split}
& \dfrac{dq^i}{dt}=\dot q^i, \\
& \dfrac{d\dot q^i}{dt}={}-\Gamma_{hk}^i\,\dot q^h\,\dot q^k.
\end{split}
\right. \label{XG}
  \end{equation}

We can consider $\mathbf X_\lambda=\mathbf X_G+\mathbf A+\mathbf
R_\lambda$ as a {\it family} of vector f\/ields, depending on the
Lagrangian multipliers. However, it is a remarkable fact that we
can obtain an explicit form of them, as functions of the
kinematical states $(q,\dot q)$ only.

\begin{theorem}
Let $[G^{ab}]$ be the symmetric matrix defined by
\begin{equation}\label{GAB}
G^{ab} \doteq g^{ij}\,C^a_i\,C^b_j.
  \end{equation}
Let  $[G_{ab}]=[G^{ab}]^{-1}$ be its inverse, and $C^{ai}\doteq
g^{ij}\,C^a_j$. If the constraints are ideal, then the Lagrangian
multipliers and the Lagrangian reactive forces are well determined
functions of $(q,\dot q)$:
  \begin{gather}\label{la}
\lambda_a(q,\dot q)= G_{ab}\left( C^b_i\,(\Gamma^i_{hk}\,\dot
q^h\,\dot q^k-A^i) -\dot q^i\,\partial_iC^b \right),
\\
\label{laR} R^i(q,\dot q)= G_{ab}\, C^{ai}\left(
C^b_j\,(\Gamma^j_{hk}\,\dot q^h\,\dot q^k-A^j) -\dot
q^j\,\partial_jC^b \right).
\end{gather}
\end{theorem}

\begin{proof}
First of all, we observe that the matrix $[G^{ab}]$ is regular,
since the vector f\/ields $\mathbf C^a$ are independent. Hence,
the inverse matrix $[G_{ab}]$ is well def\/ined. In order to
satisfy the constraints, the vector f\/ield $\mathbf X_\lambda$
must be tangent to the constraint submanifold  $C$. This condition
is expressed  by equations
    \begin{equation}
\langle \mathbf X_\lambda,dC^a\rangle =0,
    \end{equation}
to be satisf\/ied at least on $C$. In components, these equations
read
     \[
  \dot q^i\,\partial_iC^a+ (A^i-\Gamma^i_{hk}\,\dot q^h\,\dot q ^k+R^i)\,C^a_i=0,
    \]
i.e.,
  \[
R^iC_i^a=C^a_i\,(\Gamma^i_{hk}\,\dot q^h\,\dot q^k-A^i)-\dot
q^i\,\partial_iC^a.
  \]
Note that the right hand side
  \begin{equation}
\Lambda^a(q,\dot q)\doteq C^a_i\,(\Gamma^i_{hk}\,\dot q^h\,\dot
q^k-A^i)-\dot q^i\,\partial_iC^a \label{Lambda}
  \end{equation}
is a known function of $(q,\dot q)$. Then equation
$R^iC_i^a=\Lambda^a$ assumes the form
$\lambda_b\,C^a_i\,C^b_j\,g^{ij}=\Lambda^a$, i.e.,
$\lambda_b\,G^{ab}=\Lambda^a$. By applying the inverse matrix
$[G_{ab}]$ we get equation~(\ref{la}) and equation~(\ref{laR}).
\end{proof}

The explicit form (\ref{laR}) of the reactive forces allows us to
state

\begin{theorem}\label{t:D}
The actual motions of a mechanical system with regular and ideal
non-holonomic constraints represented by a submanifold $C\subset
TQ$ are the integral curves based on $C$ of the vector field
\begin{equation}\label{primoX}
\mathbf D\doteq \mathbf X_G+\mathbf A+\mathbf R,
\end{equation}
where the components $R^i$ of the vertical vector $\mathbf R$ are
defined by $(\ref{laR})$.  If we introduce the symbols
\begin{equation}\label{PIJ}
\pi^{ij}\doteq G_{ab}C^{ai}\,C^{bj}, \qquad C^{ai}\doteq
g^{ij}\,C^a_j, \qquad C_a^i\doteq G_{ab}\,C^{bi},
\end{equation}
then the explicit expressions of the first-order differential
system associated with $\mathbf D$ and of the reactive forces are
 \begin{equation}\label{DD}
 \mathbf D=
\left\{
 \begin{split}
& \dfrac{dq^i}{dt}=\dot q^i, \\
& \dfrac{d\dot q^i}{dt}=(g^{ij}-\pi^{ij})\,(A_j-\Gamma_{hkj}\,\dot
q^h\,\dot q^k)
 -\dot q^j\,\partial_jC^a\,C_a^i,
\end{split}
\right.
 \end{equation}
and
\begin{equation}\label{Rex}
R^i= \pi^{ij}\,(\Gamma_{hkj}\,\dot q^h\,\dot q^k-A_j) -\dot
q^j\,\partial_jC^a\,C_a^i,
\end{equation}
respectively.
\end{theorem}

\begin{remark}\label{r:SR}
This last theorem provides a second `\emph{recipe}' for writing
the dynamical equations for non-holonomic systems with linear or
non-linear ideal constraints:
\begin{enumerate}\itemsep=0pt
\item Choose Lagrangian coordinates $(q^i)$, write the kinetic
energy of the system $K=\tfrac 12\; g_{ij}\,\dot q^i\dot q^j$,
extract the $n\times n$ matrix $[g_{ij}]$, and compute the inverse
matrix $[g^{ij}]$. \item Take the constraint equations $C^a(q,\dot
q)=0$ and compute, in the order, the following matrices:
\begin{gather*}
[C_i^a]\doteq\left[ \dfrac{\partial C^a}{\partial \dot q^i}
\right]\qquad \hbox{(seek the singular states)},
\\
[C^{ai}]\doteq [g^{ij}\,C^a_j],
\\
[G^{ab}]=[G^{ba}]\doteq [g^{ij}\,C^a_i\,C^b_j]=[C^{ai}\,C^b_i],
\\
[G_{ab}]\doteq [G^{ab}]^{-1},
\\
[C_a^i]\doteq [G_{ab}\,C^{bi}],
\\
[\pi^{ij}]=[\pi^{ji}]\doteq[C^{ai}\,C_a^j],
\\
[g^{ij}-\pi^{ij}],
\\
[\partial_iC^a]\doteq\left[\dfrac{\partial C^a}{\partial q^i}
\right],
\\
[\partial_jC^a\,C_a^i].
\end{gather*}
\item Write the Lagrange equations for the free motions in the
form
\[
g_{ij}\,\ddot q^j=A_i-\Gamma_{hki}\,\dot q^h\,\dot q^k,
\]
and keep in evidence the functions
\[
L_i(q,\dot q)\doteq A_i-\Gamma_{hki}\,\dot q^h\,\dot q^k.
\]
\item Compute the vector
\[
D^i\doteq (g^{ij}-\pi^{ij})\,L_j-\dot q^j\,\partial_jC^a\,C_a^i.
\]
\item Write the dif\/ferential system (\ref{DD}),
\[
\mathbf D= \left\{
 \begin{split}
& \dfrac{dq^i}{dt}=\dot q^i, \\
& \dfrac{d\dot q^i}{dt}=D^i.
\end{split}
\right.
\]
\item Its solutions $q^i(t)$, $\dot q^i(t)$, with initial
conditions belonging to $C$, describe the actual motions of the
system. \item Evaluate the reactive forces along any actual motion
by means of equation~(\ref{Rex}),
\begin{equation}
R^i= {}-\pi^{ij}\,L_j -\dot q^j\,\partial_jC^a\,C_a^i.
\end{equation}
\end{enumerate}
\end{remark}

\begin{remark} \label{SR1}
{\bf The case of single constraint equation} $C(q,\dot q)=0$. In
this case the above-given recipe can be applied by setting
$a=b=1$. Items 1 and 2 of the general recipe still hold. However,
since some of the above matrices reduces to scalar functions or to
vectors, the index $1$ can be omitted or replaced by $\ast$:
\begin{gather*}
[C_i]\doteq\left[ \dfrac{\partial C}{\partial \dot q^i}
\right],
\\
[C^{i}]\doteq [g^{ij}\,C_j],
\\
G\doteq [g^{ij}\,C_i\,C_j]=C^{i}\,C_i,
\\
C_\ast^i\doteq G^{-1}\,C^i,
\\
[\pi^{ij}]\doteq[C^{i}\,C_\ast^j]=G^{-1}\,[C^i\,C^j],
\\
[g^{ij}-\pi^{ij}],
\\
\partial_iC\doteq\dfrac{\partial C}{\partial q^i},
\\
[\partial_jC\,C_\ast^i].
\end{gather*}
Then, follows items 3--7 of the general recipe.
\end{remark}

\section{Illustrative examples}

As shown above, for writing the dynamical equations of a
non-holonomic system we can apply two methods: the f\/irst method
is established by Theorem \ref{t:Z} and the corresponding recipe
is illustrated in Remark \ref{r:Z}; the second method is
established by Theorem \ref{t:D} and the corresponding recipe is
illustrated in Remark \ref{r:SR}.

Let us see how these two methods work by concrete examples. We
begin with two paradigmatic and simple examples of linear
non-holonomic constraints, the `\emph{skate}' and the
`\emph{vertical rolling disc}'. Then, we shall consider two more
demanding examples: `\emph{two co-axial rolling discs}' and
`\emph{two points with parallel velocities}'. This last one is a
genuine non-linear non-holonomic system.

In illustrating examples of application of a theory it is not
customary, in general, to provide detailed calculations -- which
usually are left to the reader. Here, however, it is worthwhile to
disregard such a custom in order to compare the ef\/fectiveness
of the two methods (mainly the length of the calculations)
applied to a same mechanical system.

\subsection{The skate}

This mechanical system is made of a homogeneous rod (material
segment) sliding without friction on a plane\footnote{Quite
similar classical examples are that of two material points linked
by a massless rigid segment \cite{Gantmacher}, p.~23 and 63, and
the \v Caplygin sleigh  \cite {Neimark}, Ch.~III, \S~3, Examples~2
\&~5, and  Ch.~V, \S~4.  Another example of this kind is examined
in~\cite{Monforte}, \S~4.1 \& \S~4.2.}. The conf\/igurations of
the skate are determined by the Cartesian coordinates $(x,y)$ of
the center of mass (i.e., of the segment) $G$ and by the angle
$\theta$ of the rod w.r.to the $x$-axis. The conf\/iguration
manifold $Q$ is $\mathbb R\times \mathbb S_1$ and natural
\emph{ordered} Lagrangian coordinates are
$(q^1,q^2,q^3)=(x,y,\theta)$. The velocity $\mathbf v_G=[\dot
x,\dot y]$ of the mass-center is constrained to be parallel to the
rod. This constraint is then represented by a single equation:
  \begin{equation}\label{SK1}
\dot x \sin\theta-\dot y\cos\theta=0.
  \end{equation}
The kinetic energy is given by $K=\tfrac 12\, m\, (\dot x^2+\dot
y^2)+ \tfrac 12 \,I\, \dot\theta$, where $m$ and $I$ are the mass
and the moment of inertia w.r.to $G$ (i.e., w.r.to the line
orthogonal to the plane through $G$), respectively.

(i) First method. Since $\dim (TQ)=6$ and $\dim(C)=3$, we need two
parameters $(z^1,z^2)$ for describing $C$. We can consider the
parametric equations
\begin{gather*}
\dot x = z^1\,\cos\theta, \qquad \dot y = z^1\,\sin\theta, \qquad
\dot\theta = z^2.
 \end{gather*}
Then, we compute the necessary matrices and vectors:
\begin{gather*}
[g_{ij}]= \left[
\begin{array}{ccc}
m & 0 & 0 \\
0 & m & 0 \\
0 & 0 & I \\
\end{array}
\right],
\\
[\psi^i]=[z^1\,\cos\theta \; , \; z^1\, \sin\theta \; , \; z^2],
\qquad [\psi^i_\alpha]= \left[
\begin{array}{ccc}
\cos\theta & \sin\theta & 0 \\
0 & 0 & 1 \\
\end{array}
\right].
\end{gather*}
This matrix has maximal rank everywhere: the constraint is
regular.
  \[
[G_{\alpha\beta}]= \left[
\begin{array} {cc}
m  & 0 \\
0 & I  \\
\end{array}
\right], \qquad [G^{\alpha\beta}]= \left[
\begin{array}{cc}
\frac 1m  & 0 \\
0 & \frac 1I  \\
\end{array}
\right].
  \]
Since,
  \[
\left[\dfrac{\partial K}{\partial \dot q^i}\right]=[m\,\dot x
\;,\; m\,\dot y  \;,\;  I\,\dot\theta  ], \qquad
\left[\dfrac{\partial K}{\partial q^i}\right]=[0\;,\;0 \;,\; 0],
 \]
the Lagrange equations for the free motions $g_{ij}\,\ddot
q^j=A_i-\Gamma_{hki}\,\dot q^h\,\dot q^k$ read
  \[
m\,\ddot x=A_1  ,   \qquad m\,\ddot y=A_2, \qquad I\,\ddot
\theta=A_3 .
  \]
Hence, $L_i=A_i(q,\dot q)$, $\bar Z_i=A_i(q,\psi)$ and
  \[
[Z_\alpha]=[\psi^i_\alpha\,\bar Z_i]= \left[
\begin{array} {c}
A_1 \,\cos\theta + A_2 \,\sin\theta \\
A_3 \\
\end{array}
\right] .
  \]
The dynamical equations are
 \[
\begin{split}
&\dfrac{dx}{dt} = z^1\,\cos\theta, \\
&\dfrac{dy}{dt}= z^1\,\sin\theta, \\
&\dfrac{d\theta}{dt}= z^2,
 \end{split}
 \qquad\quad
 \begin{split}
&\dfrac{dz^1}{dt}=
\dfrac 1m\; (A_1 \,\cos\theta + A_2 \,\sin\theta), \\
&\dfrac{dz^2}{dt}=\frac {A_3}I .
 \end{split}
   \]

(ii) Second method, for a single constraint equation -- Remark
\ref{SR1}: $C\doteq\dot x \sin\theta-\dot y\cos\theta=0$.
\begin{gather*}
[g^{ij}]= \left[
\begin{array}{ccc}
\tfrac 1m & 0 & 0 \\
0 & \tfrac 1m & 0 \\
0 & 0 & \tfrac 1I \\
\end{array}
\right],
\\
[C_i]\doteq\left[\dfrac{\partial C}{\partial \dot
q^i}\right]=\left[\sin\theta \;,\; {}-\cos\theta \;,\; 0\right],
\\
[C^i]\doteq[g^{ij}\,C_j]=\left[\dfrac{\sin\theta}m \;,\;
{}-\dfrac{\cos\theta}m \;,\; 0\right],
\\
G\doteq C^iC_i=\frac 1m,
\\
[C_*^i]\doteq G^{-1}[C^i]=\left[\sin\theta \;,\; {}-\cos\theta
\;,\; 0\right]=[C_i],
\\
[\pi^{ij}]\doteq G^{-1}[C^iC^j]=\dfrac 1m\, \left[
\begin{array}{ccc}
\sin^2\theta & {}-\sin\theta\;\cos\theta & 0 \\
{}-\sin\theta\;\cos\theta & \cos^2\theta & 0 \\
0 & 0 & 0 \\
\end{array}
\right],
\\
[g^{ij}-\pi^{ij}]=\dfrac 1m\, \left[
\begin{array}{ccc}
\cos^2\theta & \sin\theta\;\cos\theta & 0 \\
\sin\theta\;\cos\theta & \sin^2\theta & 0 \\
0 & 0 & \frac mI \\
\end{array}
\right],
\\
[\partial_iC]\doteq\left[\dfrac{\partial C}{\partial q^i}
\right]=\left[0 \;,\; 0 \;,\; \dot x\;\cos\theta+\dot
y\;\sin\theta \right],
\\
[\dot q^j\,\partial_jC\,C_\ast^i]=\dot\theta \,(\dot
x\;\cos\theta+\dot y\;\sin\theta) \,[C_\ast^i] .
\end{gather*}
Since
 \[
\left[\dfrac{\partial K}{\partial \dot q^i}\right]=[m\,\dot x\;,\;
m\,\dot y \,,\, I\,\dot\theta], \qquad \left[\dfrac{\partial
K}{\partial q^i}\right]=[0\;,\;0 \;,\; 0],
 \]
the Lagrange equations read
  \[
m\,\ddot x= A_1,\qquad m\,\ddot y= A_2, \qquad I\,\ddot\theta=A_3.
   \]
Hence, $L_i=A_i$. We have all the ingredients for computing the
vector $D^i\doteq (g^{ij}-\pi^{ij})\,L_j-\dot
q^j\,\partial_jC\,C_\ast^i$:
\begin{gather*}
D^1= \dfrac{\cos\theta}{m}\,\left(A_1\,\cos\theta+A_2\,\sin\theta
\right)- \dot\theta \,(\dot x\;\cos\theta+\dot
y\;\sin\theta)\,\sin\theta,
\\
D^2=\dfrac{\sin\theta}{m}\,\left(
A_1\,\cos\theta+A_2\,\sin\theta\right)+ \dot\theta \,(\dot
x\;\cos\theta+\dot y\;\sin\theta)\,\cos\theta,
\\
D^3=\frac 1I\,A_3,
\end{gather*}
and the dif\/ferential system (\ref{DD}) reads
 \[
\begin{split}
\dfrac{dx}{dt}=\dot x, \\
\dfrac{dy}{dt}=\dot y, \\
\dfrac{d\theta}{dt}=\dot \theta,
\end{split}
\qquad\quad
\begin{split}
&\dfrac{d\dot x}{dt}= \dfrac{\cos\theta}{m}\,\left(A_1\,\cos\theta+A_2\,\sin\theta \right)- \dot\theta \,(\dot x\;\cos\theta+\dot y\;\sin\theta)\,\sin\theta, \\
&\dfrac{d\dot y}{dt}=\dfrac{\sin\theta}{m}\,\left(
A_1\,\cos\theta+A_2\,\sin\theta\right)+
\dot\theta \,(\dot x\;\cos\theta+\dot y\;\sin\theta)\,\cos\theta,\\
&\dfrac{d\dot \theta}{dt}=I^{-1}\,A_3.
\end{split}
\]


\subsection{The vertical rolling disc}

A material disc of radius $R$ running on a plane is kept
perpendicular to it by  massless and frictionless devices. The
conf\/iguration manifold is $Q_4=\mathbb R^2\times \mathbb
S_1\times\mathbb S_1$, with coordinates
$(q^1,q^2,q^3,q^4)=(x,y,\theta,\psi)$, where $(x,y)$ are Cartesian
coordinates of the center $P$ of the disc (i.e., of the point $C$
in contact with the plane), $\theta$ a rotation angle of the disk
around its axis, and $\psi$ an angle giving the orientation of the
axis (see Figure 4, with $\theta=\theta_1$). Constraint: the disc
rolls on the plane without sliding. Let $(\mathbf i, \mathbf j,
\mathbf k)$ be the unitary vectors associated with the
$(x,y,z)$-axes. The unitary vector $\mathbf k$ is associated with
the oriented angle $\psi$. Let $\mathbf u$ be the unitary vector
associated with the oriented angle $\theta$. Then, $\mathbf
u=\cos\psi\;\mathbf j-\sin\psi\;\mathbf i$. The angular velocity
$\mathbf \omega$ is given by $\omega=\dot\theta\,\mathbf u+
\dot\psi\,\mathbf k$. The velocity $\mathbf v_C$ of the point $C$
is given by $\mathbf v_C=\mathbf v_P+\omega\times PC$, where
$PC={}-R\,\mathbf k$. Hence,
\[
\mathbf v_C=\dot x\,\mathbf i+\dot y\,\mathbf
j-(\dot\theta\,\mathbf u+ \dot\psi\,\mathbf k)\times R\,\mathbf
k=\dot x\,\mathbf i+\dot y\,\mathbf j-R\,\dot\theta\,\mathbf
u\times\mathbf k.
\]
Since $\mathbf u\times \mathbf k=\cos\psi\,\mathbf j\times\mathbf
k-\sin\psi\,\mathbf i\times \mathbf k$, we get
\begin{gather*}
\mathbf v_C=\dot x\,\mathbf i+\dot y\,\mathbf j-R\,\dot\theta\,(\cos\psi\,\mathbf j\times\mathbf k-\sin\psi\,\mathbf i\times \mathbf k)=\dot x\,\mathbf i+\dot y\,\mathbf j-R\,\dot\theta\,(\cos\psi\,\mathbf i+\sin\psi\,\mathbf j)\\
\phantom{\mathbf v_C} {}=(\dot x-
R\,\dot\theta\,\cos\psi)\,\mathbf i+(\dot
y-R\,\dot\theta\,\sin\psi)\,\mathbf j.
\end{gather*}
The kinematical constraint $\mathbf v_C=0$  is then represented by
the following two linear equations
 \begin{gather}
 C^1\doteq\dot x- R \,\cos\psi\,\dot\theta=0, \nonumber\\
 C^2\doteq\dot y - R \, \sin\psi  \,\dot\theta=0.
\label{RDC}
\end{gather}

(i) First method. Assume that the center of mass of the disc
coincides with its geometrical center. Then the kinetic energy is
given by
  \begin{equation}
K={\tfrac 12} \, m (\dot x^2+\dot y ^2)+ {\tfrac 12}\,  (A \dot
\theta^2+B\dot\psi^2), \label{Kk}
  \end{equation}
where $m$ is the mass, $A$ and $B$ are the moments of inertia
w.r.to the axis of rotation and a~diameter, respectively. Thus,
  \[
[g_{ij}] = \left[
\begin{array} {cccc}
m & 0 & 0 & 0 \\
0 & m & 0 & 0 \\
0 & 0 & A & 0 \\
0 & 0 & 0 & B
\end{array}
\right].
  \]
From the constraint equations (\ref{RDC}) we get the parametric
equations
 \begin{equation} 
 \begin{split}
& \dot x=R\,\cos\psi\,z^1, \\
& \dot y=R\,\sin\psi\,z^1, \\
\end{split}
\qquad\quad
 \begin{split}
& \dot\theta=z^1, \\
& \dot\psi=z^2.
\end{split}
\label{RDCP}
\end{equation}
Thus,
\begin{gather*}
[\psi^i]=[R\,\cos\psi\,z^1\;,\; R\,\sin\psi\,z^1 \;,\; z^1\; ,\;
z^2  ],
\\
[\psi^i_\alpha]= \left[
\begin{array}{cccc}
R\,\cos\psi & R\,\sin\psi & 1 & 0 \\
0 & 0 & 0 & 1 \\
\end{array}
\right] \qquad (\hbox{$\alpha=1,2$, index of line}).
\end{gather*}
This matrix has maximal rank, thus the constraint is regular. It
follows that
\begin{gather*}
G_{\alpha\beta}=g_{ij}\psi^i_\alpha\psi^j_\beta  =
m\,\psi^1_\alpha\psi^1_\beta +m\,\psi^2_\alpha\psi^2_\beta
+A\,\psi^3_\alpha\psi^3_\beta +B\,\psi^4_\alpha\psi^4_\beta,
\\
G_{11}=m\,\psi^1_1\psi^1_1 +m\,\psi^2_1\psi^2_1
+A\,\psi^3_1\psi^3_1 +B\,\psi^4_1\psi^4_1 = mR^2\cos^2\psi
+mR^2\sin^2\psi +A=mR^2\!+A,\!
\\
G_{12}=m\,\psi^1_1\psi^1_2 +m\,\psi^2_1\psi^2_2
+A\,\psi^3_1\psi^3_2 +B\,\psi^4_1\psi^4_2 =0,
\\
G_{22}=m\,\psi^1_2\psi^1_2 +m\,\psi^2_2\psi^2_2
+A\,\psi^3_2\psi^3_2 +B\,\psi^4_2\psi^4_2=B,
\\
[G_{\alpha\beta}]= \left[
\begin{array} {cc}
mR^2 + A & 0 \\
0 & B  \\
\end{array}
\right],
\\
[G^{\alpha\beta}]=\, \left[
\begin{array} {cc}
\dfrac 1{mR^2+A} & 0 \\
0 & \dfrac 1B  \\
\end{array}
\right].
\end{gather*}
Since,
  \[
\left[\dfrac{\partial K}{\partial \dot q^i}\right]=[m\,\dot x
\;,\; m\,\dot y  \;,\;  A\,\dot\theta \;,\; B\,\dot\psi ], \qquad
\left[\dfrac{\partial K}{\partial q^i}\right]=[0\;,\;0 \;,\; 0
\;,\; 0],
 \]
the Lagrange equations for the free motions $g_{ij}\,\ddot
q^j=A_i-\Gamma_{hki}\,\dot q^h\,\dot q^k$ read
  \[
m\,\ddot x=A_1  , \qquad m\,\ddot y=A_2, \qquad A\,\ddot
\theta=A_3 , \qquad B\,\ddot \psi=A_4 .
  \]
This shows that $L_i=A_i$. Hence, $\bar Z_i=A_i(q,\psi)$, and
\begin{gather*}
[Z_\alpha]= [\bar Z_i\,\psi^i_\alpha] = \left[
\begin{array}{c}
A_1\,R\,\cos\psi +A_2\,R\,\sin\psi + A_3  \\
A_4
\end{array}
\right],
\\
[Z^\alpha]=[G^{\alpha\beta}\,Z_\beta] = \left[
\begin{array}{c}
\dfrac{A_1\,R\,\cos\psi +A_2\,R\,\sin\psi + A_3}{mR^2+A} \\[12pt]
\dfrac{A_4}B
\end{array}
\right].
\end{gather*}
Thus, the dynamical equations are
\begin{equation}
\begin{split}
&\dfrac{dx}{dt} =R\,\cos\psi\,z^1, \\
&\dfrac{dy}{dt} =R\,\sin\psi\,z^1, \\
&\dfrac{d\theta}{dt} =z^1, \\
&\dfrac{d\psi}{dt} =z^2,
\end{split}
\qquad\quad
\begin{split}
&\dfrac{dz^1}{dt} = (mR^2+A)^{-1}\,(A_1\,R\,\cos\psi +A_2\,R\,\sin\psi + A_3), \\
&\dfrac{dz^2}{dt}=B^{-1}\,A_4.
\end{split}
\label{RDZ1}
 \end{equation}

(ii) Second method. Recall the constraint equations (\ref{RDC}).
Then,
\begin{gather*}
[C^a_i]\doteq\left[ \dfrac{\partial C^a}{\partial \dot q^i}
\right]= \left[
\begin{array} {cccc}
1 & 0 & {}-R\,\cos\psi & 0 \\[8pt]
0 & 1 & {}-R\,\sin\psi & 0
\end{array}
\right] ,
\\
[C^{ai}]\doteq \left[g^{ij} \,C^a_j \right]= \left[
\begin{array} {cccc}
\frac 1m & 0 & {}-\frac RA \,\cos\psi & 0 \\[8pt]
0 & \frac 1m & {}-\frac RA \,\sin\psi & 0
\end{array}
\right] ,
\\
[G^{ab}]\doteq\left[C^{ai}C^b_i \right]= \left[
\begin{array} {cc}
\frac 1m +\frac{R^2}A\,\cos^2\psi & \frac{R^2}{A}\,\sin\psi\cos\psi \\[8pt]
\frac{R^2}{A}\,\sin\psi\cos\psi &  \frac 1m
+\frac{R^2}A\,\sin^2\psi
\end{array}
\right] ,
\\
G\doteq \det[G^{ab}]=\frac1{m^2}+\frac{R^2}{mA},\qquad
G^{-1}=\dfrac{m^2A}{mR^2+A},
\\
[G_{ab}]= G^{-1} \left[
\begin{array} {cc}
\frac 1m +\frac{R^2}A\,\sin^2\psi & {}-\frac{R^2}{A}\,\sin\psi\cos\psi \\[8pt]
{}-\frac{R^2}{A}\,\sin\psi\cos\psi &  \frac 1m
+\frac{R^2}A\,\cos^2\psi
\end{array}
\right],
\\
[C_a^i]= [G_{ab}\,C^{bi}] =G^{-1} \left[
\begin{array} {cccc}
\frac 1m\,(\frac 1m +\frac{R^2}A\,\sin^2\psi)& {}-\frac{R^2}{mA}\,\sin\psi\cos\psi  & *C_1^3 & 0 \\
{}-\frac{R^2}{mA}\,\sin\psi\cos\psi & \frac 1m\,(\frac 1m
+\frac{R^2}A\,\cos^2\psi) & *C_2^3 & 0
\end{array}
\right],
\end{gather*}
where
\begin{gather*}
*C_1^3=\left(\frac 1m +\frac{R^2}A\,\sin^2\psi\right)\left(-\frac
RA \,\cos\psi \right) +
\left(-\frac{R^2}{A}\,\sin\psi\cos\psi\right)\left(-\frac RA
\,\sin\psi \right)
\\
\phantom{*C_1^3}{}=-\frac R{mA}
\,\cos\psi-\frac{R^3}{A^2}\,\sin^2\psi\,\cos\psi +
\frac{R^3}{A^2}\,\sin^2\psi\cos\psi={}-\frac R{mA} \,\cos\psi,
\\
*C_2^3= \left(-\frac{R^2}{A}\,\sin\psi\cos\psi \right)\left(-\frac
RA \,\cos\psi\right)+
\left( \frac 1m +\frac{R^2}A\,\cos^2\psi\right)\left(-\frac RA \,\sin\psi\right) \\
\phantom{*C_2^3}{}=\left(-\frac{R^2}{A}\,\sin\psi\cos\psi
\right)\left(-\frac RA \,\cos\psi\right)+
\left( \frac 1m +\frac{R^2}A\,\cos^2\psi\right)\left(-\frac RA \,\sin\psi\right) \\
\phantom{*C_2^3}{}=-\frac R{mA}\,\sin\psi.
\end{gather*}
Hence,
\[
[C_a^i]= G^{-1} \left[
\begin{array} {cccc}
\frac 1m\,(\frac 1m +\frac{R^2}A\,\sin^2\psi)& {}-\frac{R^2}{mA}\,\sin\psi\cos\psi  & {}-\frac R{mA} \,\cos\psi & 0 \\[12pt]
{}-\frac{R^2}{mA}\,\sin\psi\cos\psi & \frac 1m\,(\frac 1m
+\frac{R^2}A\,\cos^2\psi) & {}-\frac R{mA}\,\sin\psi & 0
\end{array}
\right].
\]
Let us compute $\pi^{ij}\doteq C^{ai}\,C_a^j=\pi^{ji}$. Let us set
$*\pi^{ij}\doteq G\,\pi^{ij}$. Then,
\[
[\pi^{ij}]= G^{-1} \left[
\begin{array} {cccc}
\frac 1{m^2}\,(\frac 1m +\frac{R^2}A\,\sin^2\psi) &
{}-\frac{R^2}{m^2A}\,\sin\psi\cos\psi & *\pi^{14}
& 0 \\[12pt]
{}-\frac{R^2}{m^2A}\,\sin\psi\cos\psi & \frac 1{m^2}\,(\frac 1m
+\frac{R^2}A\,\cos^2\psi) & *\pi^{24}
& 0  \\[12pt]
 *\pi^{31} &  *\pi^{32} &  *\pi^{33} &  0 \\[12pt]
 *\pi^{41} & *\pi^{42} & 0 & 0
\end{array}
\right],
\]
where
\begin{gather*}
*\pi^{13} =\frac 1m\left(\frac 1m
+\frac{R^2}A\,\sin^2\psi\right)\left(-\frac RA\,\cos\psi\right) +
\left(-\frac{R^2}{mA}\,\sin\psi\cos\psi\right)\left(-\frac RA\,\sin\psi\right) \\
\phantom{*\pi^{13}}{} =\frac 1m\left(\frac 1m
+\frac{R^2}A\,\sin^2\psi\right)\left(-\frac RA\,\cos\psi\right) +
\frac{R^3}{mA^2}\,\sin^2\psi\cos\psi =-\frac{R}{m^2A}\,\cos\psi,
\\
*\pi^{23}
=\left(-\frac{R^2}{mA}\,\sin\psi\cos\psi\right)\left(-\frac
RA\,\cos\psi\right) +
\frac 1m\left(\frac 1m +\frac{R^2}A\,\cos^2\psi\right)\left(-\frac RA\,\sin\psi\right) \\
\phantom{*\pi^{23}}{} =-\frac{R}{m^2A}\,\sin\psi.
\end{gather*}
It follows that
\begin{gather*}
[\pi^{ij}]= G^{-1} \left[
\begin{array} {cccc}
\frac 1{m^2}\,(\frac 1m +\frac{R^2}A\,\sin^2\psi) &
{}-\frac{R^2}{m^2A}\,\sin\psi\cos\psi &
{}-\frac{R}{m^2A}\,\cos\psi
& 0 \\[12pt]
{}-\frac{R^2}{m^2A}\,\sin\psi\cos\psi & \frac 1{m^2}\,(\frac 1m
+\frac{R^2}A\,\cos^2\psi) & {}-\frac{R}{m^2A}\,\sin\psi
& 0  \\[12pt]
 {}-\frac{R}{m^2A}\,\cos\psi  & {}-\frac{R}{m^2A}\,\sin\psi & \frac{R^2}{mA^2} & 0 \\[12pt]
 0 & 0 & 0 & 0
\end{array}
\right],
\\
[g^{ij}-\pi^{ij}]= G^{-1} \left[
\begin{array} {cccc}
\frac Gm-\frac 1{m^2}\,(\frac 1m +\frac{R^2}A\,\sin^2\psi) &
\frac{R^2}{m^2A}\,\sin\psi\cos\psi & \frac{R}{m^2A}\,\cos\psi
& 0 \\[12pt]
\frac{R^2}{m^2A}\,\sin\psi\cos\psi & \frac Gm-\frac 1{m^2}\,(\frac
1m +\frac{R^2}A\,\cos^2\psi) & \frac{R}{m^2A}\,\sin\psi
& 0  \\[12pt]
\frac{R}{m^2A}\,\cos\psi  & \frac{R}{m^2A}\,\sin\psi &
 \frac GA-\frac{R^2}{mA^2} & 0 \\[12pt]
 0 & 0 & 0 & \frac GB
\end{array}
\right]
\\
= \dfrac{m^2A}{mR^2+A} \left[
\begin{array} {cccc}
\frac Gm-\frac 1{m^2}\,(\frac 1m +\frac{R^2}A\,\sin^2\psi) &
\frac{R^2}{m^2A}\,\sin\psi\cos\psi & \frac{R}{m^2A}\,\cos\psi
& 0 \\[12pt]
\frac{R^2}{m^2A}\,\sin\psi\cos\psi & \frac Gm-\frac 1{m^2}\,(\frac
1m +\frac{R^2}A\,\cos^2\psi) & \frac{R}{m^2A}\,\sin\psi
& 0  \\[12pt]
\frac{R}{m^2A}\,\cos\psi  & \frac{R}{m^2A}\,\sin\psi &
 \frac GA-\frac{R^2}{mA^2} & 0 \\[12pt]
 0 & 0 & 0 & \frac GB
\end{array}
\right]
\\
= \dfrac{1}{mR^2+A} \left[
\begin{array} {cccc}
R^2\,\cos^2\psi & R^2\,\sin\psi\cos\psi & R\,\cos\psi
& 0 \\[12pt]
R^2\,\sin\psi\cos\psi & R^2\,\sin^2\psi & R\,\sin\psi
& 0  \\[12pt]
R\,\cos\psi  & R\,\sin\psi &
1 & 0 \\[12pt]
 0 & 0 & 0 & \frac {mR^2+A}B
\end{array}
\right] .
\end{gather*}
Recall once more equations (\ref{RDC}). Then,
\[
[\partial_iC^a]\doteq \left[\dfrac{\partial C^a}{\partial
q^i}\right] \left[
\begin{array} {cccc}
0 & 0 & 0 &  R \,\sin\psi\,\dot\theta \\[8pt]
0 & 0 & 0 &  {}-R \,\cos\psi\,\dot\theta
\end{array}
\right] , \qquad [\dot q^i\,\partial_iC^a]= \left[
\begin{array} {c}
  R \,\sin\psi\,\dot\theta\,\dot\psi \\[8pt]
 {}-R \,\cos\psi\,\dot\theta\,\dot\psi
\end{array}
\right].
\]
Let us set $X^i\doteq\dot q^j\,\partial_jC^a\,C_a^i$. Then,
\begin{gather*}
GX^1=
 (R \,\sin\psi\,\dot\theta\,\dot\psi)\left(\frac 1m\left(\frac 1m +\frac{R^2}A\,\sin^2\psi\right)\right)
 +
(-R
\,\cos\psi\,\dot\theta\,\dot\psi)\left(-\frac{R^2}{mA}\,\sin\psi\cos\psi\right)
 \\
\phantom{GX^1}{}=\frac R{m^2} \,\sin\psi\,\dot\theta\,\dot\psi
+\frac {R^3}{mA} \,\sin\psi\,\dot\theta\,\dot\psi =\frac R{m^2}
\,\sin\psi\,\dot\theta\,\dot\psi\,\left(1+ \frac{mR^2}{A}  \right)
\\
\phantom{GX^1}{}=\frac {R \,(mR^2+A) }{m^2A}
\,\sin\psi\,\dot\theta\,\dot\psi,
\\
GX^2=
 (R \,\sin\psi\,\dot\theta\,\dot\psi)\left(-\frac{R^2}{mA}\,\sin\psi\cos\psi\right)
 +
 (-R \,\cos\psi\,\dot\theta\,\dot\psi)\left(\frac 1m\left(\frac 1m +\frac{R^2}A\,\cos^2\psi\right)\right)
 \\
\phantom{GX^2}{}= -\frac R{m^2} \,\cos\psi\,\dot\theta\,\dot\psi
-\frac {R^3}{mA} \,\cos\psi\,\dot\theta\,\dot\psi =-\frac {R
\,(mR^2+A) }{m^2A} \,\cos\psi\,\dot\theta\,\dot\psi,
\\
GX^3=
 (R \,\sin\psi\,\dot\theta\,\dot\psi)\left(-\frac R{mA} \,\cos\psi\right)
 +
 (-R \,\cos\psi\,\dot\theta\,\dot\psi)\left(-\frac R{mA} \,\sin\psi\right)=0,
\\
GX^4=
 (R \,\sin\psi\,\dot\theta\,\dot\psi)(0)
 +
 (-R \,\cos\psi\,\dot\theta\,\dot\psi)(0)=0.
\end{gather*}
Since $G^{-1}=\frac{m^2A}{mR^2+A}$, we get
\begin{gather*}
[X^i]=[\dot q^j\,\partial_jC^a\,C_a^i]=G^{-1\,} \left[ \frac {R
\,(mR^2+A) }{m^2A} \,\sin\psi\,\dot\theta\,\dot\psi \;,\; -\frac
{R \,(mR^2+A) }{m^2A} \,\cos\psi\,\dot\theta\,\dot\psi \;,\; 0
\;,\; 0\right]
\\
\phantom{[X^i]}{}=R\,\dot\theta\,\dot\psi \left[\sin\psi \;,\;
-\cos\psi\;,\; 0 \;,\; 0\right].
\end{gather*}
Now we are able to compute the components $D^i\doteq
(g^{ij}-\pi^{ij})\,L_j-\dot q^j\,\partial_jC^a\,C_a^i$ of the
vector $\mathbf D$:
\begin{gather*}
D^1= \dfrac{1}{mR^2+A} \left(
A_1\,R^2\,\cos^2\psi+A_2\,R^2\,\sin\psi\cos\psi  +A_3\,R\,\cos\psi
\right) -R\,\dot\theta\,\dot\psi\,\sin\psi
\\
\phantom{D^1}{}= \dfrac{R\,\cos\psi}{mR^2+A} \left(
A_1\,R\,\cos\psi+A_2\,R\,\sin\psi  +A_3 \right)
-R\,\dot\theta\,\dot\psi\,\sin\psi ,
\\
D^2= \dfrac{1}{mR^2+A} \left( A_1\,R^2\,\sin\psi\cos\psi +
A_2\,R^2\,\sin^2\psi +A_3\,R\,\sin\psi \right)
+R\,\dot\theta\,\dot\psi\,\cos\psi
\\
\phantom{D^2}{}= \dfrac{R\,\sin\psi}{mR^2+A} \left(
A_1\,R\,\cos\psi + A_2\,R\,\sin\psi +A_3 \right)
+R\,\dot\theta\,\dot\psi\,\cos\psi,
\\
D^3=\dfrac{1}{mR^2+A} \left( A_1\,R\,\cos\psi  +A_2\, R\,\sin\psi
+A_3 \right)
\\
\phantom{D^3}{}=\dfrac{1}{mR^2+A} \left( A_1\,R\,\cos\psi  +A_2\,
R\,\sin\psi +A_3 \right),
\\
D^4=\frac{A_4}B.
\end{gather*}
The resulting dynamical system (\ref{DD}) is
\begin{equation}\label{RDD1}
\begin{split}
&\dfrac{dx}{dt}=\dot x,
\\
&\dfrac{dy}{dt}=\dot y,
\\
&\dfrac{d\theta}{dt}= \dot\theta,
\\
&\dfrac{d\psi}{dt}=\dot\psi,
\end{split}
\qquad\quad
\begin{split}
&\dfrac{d\dot x}{dt}=\dfrac{R\,\cos\psi}{mR^2+A} \left(
A_1\,R\,\cos\psi+A_2\,R\,\sin\psi  +A_3 \right)
-R\,\dot\theta\,\dot\psi\,\sin\psi,
\\
&\dfrac{d\dot y}{dt}=\dfrac{R\,\sin\psi}{mR^2+A} \left(
A_1\,R\,\cos\psi + A_2\,R\,\sin\psi +A_3 \right)
+R\,\dot\theta\,\dot\psi\,\cos\psi,
\\
&\dfrac{d\dot \theta}{dt}=\dfrac{1}{mR^2+A} \left(
A_1\,R\,\cos\psi  +A_2\, R\,\sin\psi +A_3 \right),
\\
&\dfrac{d\dot \psi}{dt}=\frac{A_4}B.
\\
\end{split}
\end{equation}
By introducing the new variables
\[
X\doteq \dfrac{mR^2+A}{R}\,x, \qquad Y\doteq \dfrac{mR^2+A}{R}\,y,
\qquad \Theta \doteq (mR^2+A)\,\theta,
\]
it assumes the more compact form
\begin{equation}\label{RDD2}
\begin{split}
&\dfrac{dX}{dt}=\dot X,
\\
&\dfrac{dY}{dt}=\dot Y,
\\
&\dfrac{d\Theta}{dt}= \dot\Theta,
\\
&\dfrac{d\psi}{dt}=\dot\psi,
\end{split}
\qquad\quad
\begin{split}
&\dfrac{d\dot X}{dt}=\cos\psi \left(
A_1\,R\,\cos\psi+A_2\,R\,\sin\psi  +A_3 \right)
-\dot\Theta\,\dot\psi\,\sin\psi,
\\
&\dfrac{d\dot Y}{dt}=\sin\psi \left( A_1\,R\,\cos\psi +
A_2\,R\,\sin\psi +A_3 \right) +\dot\Theta\,\dot\psi\,\cos\psi,
\\
&\dfrac{d\dot \Theta}{dt}= \left( A_1\,R\,\cos\psi  +A_2\,
R\,\sin\psi +A_3 \right),
\\
&\dfrac{d\dot \psi}{dt}=\frac{A_4}B.
\\
\end{split}
\end{equation}
Note that in these new variables, by considering also
\[
Z^1 \doteq (mR^2+A)\,z^1,
\]
the dif\/ferential system  (\ref{RDZ1}) obtained by the f\/irst
method reads
\begin{equation}\label{RDZ2}
\begin{split}
&\dfrac{dX}{dt} =\cos\psi\,Z^1, \\
&\dfrac{dY}{dt} =\sin\psi\,Z^1, \\
&\dfrac{d\Theta}{dt} =Z^1, \\
&\dfrac{d\psi}{dt} =z^2,
\end{split}
\qquad\quad
\begin{split}
&\dfrac{dZ^1}{dt} = A_1\,R\,\cos\psi +A_2\,R\,\sin\psi + A_3, \\
&\dfrac{dz^2}{dt}=B^{-1}\,A_4.
\end{split}
  \end{equation}
The two systems (\ref{RDD2}) and (\ref{RDZ2}) are in perfect
agreement.

The above detailed calculations show that for the rolling disc the
f\/irst method is much shorter than the second one.

\subsection{Two co-axial rolling discs}

Two identical material discs of radius $R$ running on a plain are
joined by a massless common axis, along with they can slide
without friction. The conf\/iguration manifold is $Q_6=\mathbb
R^2\times \mathbb S_1\times \mathbb S_1\times \mathbb S_1\times
\mathbb R$, with Lagrangian coordinates
$(q^1,q^2,q^3,q^4,q^5,q^6)=(x,y,\theta_1,\theta_2,\psi,a)$, where
$(x,y)$ are Cartesian coordinates of the center $P_1$ of one of
the two discs, $\theta_1$ and $\theta_2$ are the angles of
rotations around the common axis, $\psi$ is the angle giving the
orientation of the axis, and $a$ is the distance between the
centers (see Fig.~4).

\begin{figure}[th]
\begin{center}
\includegraphics[width=8cm]{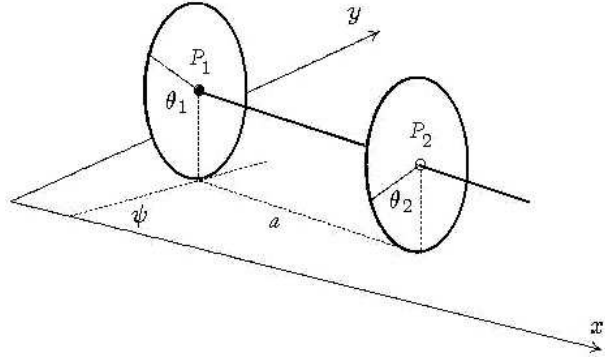}
\caption{Co-axial rolling discs.}
\end{center}
\label{Coaxial}
\end{figure}

Constraint: the discs roll on the plane without sliding. For each
disc this constraint is represented by linear equations of the
kind (\ref{RDC}),
 \begin{equation}\label{2RD}
 \begin{split}
&\dot x_1- R \,\cos\psi\,\dot\theta_1=0, \\
& \dot y_1 - R \, \sin\psi  \,\dot\theta_1=0,
\end{split}
\qquad\quad
 \begin{split}
&\dot x_2- R \,\cos\psi\,\dot\theta_2=0, \\
& \dot y_2 - R \, \sin\psi  \,\dot\theta_2=0.
\end{split}
 \end{equation}
However, the coordinates of the two centers are related by
equations
  \[
x_2= x_1+ a \;\sin\psi, \qquad y_2 = y_1-a \;\cos\psi.
 \]
By dif\/ferentiating these equations we get the link between the
velocities,
  \begin{equation}\label{RV}
\dot x_2=\dot x_1+ a \;\cos\psi\;\dot\psi+ \sin\psi\; \dot a,
\qquad \dot y_2 = \dot y_1+a \;\sin\psi\; \dot\psi-\cos\psi\;\dot
a.
 \end{equation}
By inserting these relations into equations (\ref{2RD}), with
$x_1=x$ and $y_1=y$, we get the f\/inal constraint equations
 \begin{equation}\label{2RDD}
 \begin{split}
&\dot x- R \,\cos\psi\,\dot\theta_1=0, \\
& \dot y - R \, \sin\psi  \,\dot\theta_1=0, \\
&\dot x+ a \;\cos\psi\;\dot\psi+ \sin\psi\; \dot a- R \,\cos\psi\;\dot\theta_2=0, \\
& \dot y+a \;\sin\psi\; \dot\psi-\cos\psi\;\dot a - R \, \sin\psi
\;\dot\theta_2=0.
\end{split}
 \end{equation}

Since for a single rolling disc the f\/irst method is faster, we
limit ourselves to apply  the f\/irst method to the case of two
discs.

Equations (\ref{2RDD}) show that the constraint submanifold
$C\subset TQ_6$ has dimension $12-4=8$. Hence, we need two
parameters $z^\alpha$ for the parametric representation. Let us
choose $z^1=\dot\theta_1$ and $z^2=\dot\theta_2$. This means to
solve the linear system (\ref{2RDD}) w.r.to $(\dot x, \dot y,
\dot\psi,\dot a)$. The result is
  \[
\begin{split}
&\dot x = R\;\cos\psi\;z^1, \\
&\dot y= R\;\sin\psi\;z^1, \\
&\dot\psi=\frac Ra (z^2-z^1),
\end{split}
\qquad\quad
\begin{split}
&\dot a=0, \\
&\dot \theta_1=z^1,\\
&\dot\theta_2=z^2.
\end{split}
  \]

Equation $\dot a =0$, $a=\hbox{constant}$, exhibits  the intuitive
fact that, under the pure-rolling condition, the distance $a$
between the two discs remains constant. Hence, we can reduce the
conf\/iguration manifold $Q_6$ to $Q_5= \mathbb R^2\times \mathbb
S_1\times \mathbb S_1\times \mathbb S_1$, with coordinates
$(q^1,q^2,q^3,q^4,q^5)=(x,y,\theta_1,\theta_2,\psi)$. The last
equations reduce to
\begin{equation}\label{2DPE}
\begin{split}
&\dot x = R\;\cos\psi\;z^1, \\
&\dot y= R\;\sin\psi\;z^1, \\
&\dot\psi=\frac Ra (z^2-z^1),
\end{split}
\qquad\quad
\begin{split}
&\dot \theta_1=z^1,\\
&\dot\theta_2=z^2,
\end{split}
  \end{equation}
with $a=\hbox{constant}$. Then,
  \[
[\psi^i]=\left[R\,\cos\psi\,z^1\;,\; R\,\sin\psi\,z^1 \;,\; \frac
Ra (z^2-z^1) \;,\; z^1\; ,\; z^2  \right]
   \]
and
\begin{equation}\label{2DPSI}
[\psi^i_\alpha]\doteq \left[ \dfrac{\partial \psi}{\partial
z^\alpha}  \right]= \left[
\begin{array}{ccccc}
R\,\cos\psi & R\,\sin\psi & {}-\tfrac Ra & \; 1 \; & 0 \\[8pt]
0 & 0 & \frac Ra & 0 & \; 1 \; \\
\end{array}
\right].
  \end{equation}
This matrix has maximal rank, thus the constraint submanifold $C$
is regular. The kinetic energy of the system is the sum of the
kinetic energies of the two discs. According to
equations~(\ref{Kk}) and~(\ref{RV}),
\begin{gather*}
K={\tfrac 12} \, m (\dot x_1^2+\dot y_1 ^2)+ {\tfrac 12}\, (A \dot
\theta_1^2+B\dot\psi^2) +{\tfrac 12} \, m (\dot x_2^2+\dot y_2
^2)+ {\tfrac 12}\, (A \dot \theta_2^2+B\dot\psi^2)
\\
\phantom{K}={\tfrac 12} \, m (\dot x_1^2+\dot y_1 ^2+\dot
x_2^2+\dot y_2 ^2)+ {\tfrac 12}\, (A \dot \theta_1^2+B\dot\psi^2)
+ {\tfrac 12}\, (A \dot \theta_2^2+B\dot\psi^2)\\
\phantom{K}={\tfrac 12} \, m \left(\dot x_1^2+\dot y_1 ^2+(\dot
x_1+ a \;\cos\psi\;\dot\psi)^2+(\dot y_1+a \;\sin\psi\;
\dot\psi)^2\right)+ {\tfrac A2}\, (\dot \theta_1^2+\dot
\theta_2^2)
+ B\dot\psi^2\\
\phantom{K}= {\tfrac 12} \, m \left(\dot x_1^2+\dot y_1 ^2 +\dot
x_1^2+ a^2 \;\cos^2\psi\;\dot\psi^2+ 2a\,\dot
x_1\;\cos\psi\;\dot\psi
+\dot y_1^2 + a^2 \;\sin^2\psi\; \dot\psi^2\right.  \\
\left.\phantom{K=}+2a\,\dot y_1\;\sin\psi\; \dot\psi \right)+
{\tfrac A2}\, (\dot \theta_1^2+\dot \theta_2^2)
+ B\dot\psi^2 \\
\phantom{K}={\tfrac 12} \, m \left(2\dot x_1^2+2\dot y_1^2 + a^2
\;\dot\psi^2 +2a\,\dot\psi\;(\dot x_1\;\cos\psi+\dot
y_1\;\sin\psi) \right) + \tfrac A2\, (\dot \theta_1^2+\dot
\theta_2^2)
+ B\dot\psi^2 \\
\phantom{K}= m \,(\dot x^2+\dot y^2)+ \tfrac A2\,(\dot
\theta_1^2+\dot \theta_2^2)+ (\tfrac
12\,m\,a^2+B)\,\dot\psi^2+m\,a\,\dot\psi\;(\dot x\;\cos\psi+\dot
y\;\sin\psi).
\end{gather*}
Thus,
\[
[g_{ij}]= \left[
\begin{array}{ccccc}
2m & 0 & 0 & 0 & \frac{ma}2\,\cos\psi \\
0 & 2m & 0 & 0 & \frac{ma}2\,\sin\psi \\
0 & 0 & A & 0 & 0 \\
0 & 0 & 0 & A & 0 \\
\frac{ma}2\,\cos\psi & \frac{ma}2\,\sin\psi & 0 & 0 & ma^2 +2B \\
\end{array}
\right].
\]
Moreover,
\begin{gather*}
G_{\alpha\beta}=g_{ij}\,\psi^i_\alpha\,\psi^j_\beta = g_{11}
\,\psi^1_\alpha\,\psi^1_\beta +g_{22}
\,\psi^2_\alpha\,\psi^2_\beta +g_{33}
\,\psi^3_\alpha\,\psi^3_\beta +g_{44}
\,\psi^4_\alpha\,\psi^4_\beta +g_{55}
\,\psi^5_\alpha\,\psi^5_\beta
\\
\phantom{G_{\alpha\beta}=}{} +2g_{12}\,\psi^1_\alpha\,\psi^2_\beta
+2g_{13}\,\psi^1_\alpha\,\psi^3_\beta
+2g_{14}\,\psi^1_\alpha\,\psi^4_\beta
+2g_{15}\,\psi^1_\alpha\,\psi^5_\beta
+2g_{23}\,\psi^2_\alpha\,\psi^3_\beta
+2g_{24}\,\psi^2_\alpha\,\psi^4_\beta\\
\phantom{G_{\alpha\beta}=}{} +2g_{25}\,\psi^2_\alpha\,\psi^5_\beta
+2g_{34}\,\psi^3_\alpha\,\psi^4_\beta
+2g_{35}\,\psi^3_\alpha\,\psi^5_\beta
+2g_{45}\,\psi^4_\alpha\,\psi^5_\beta\\
\phantom{G_{\alpha\beta}}{} = 2m \,\psi^1_\alpha\,\psi^1_\beta +2m
\,\psi^2_\alpha\,\psi^2_\beta +A \,\psi^3_\alpha\,\psi^3_\beta +A
\,\psi^4_\alpha\,\psi^4_\beta
+(ma^2+2B) \,\psi^5_\alpha\,\psi^5_\beta\\
\phantom{G_{\alpha\beta}=}{}
+ma\,\cos\psi\,\psi^1_\alpha\,\psi^5_\beta
+ma\,\sin\psi\,\psi^2_\alpha\,\psi^5_\beta,
\\
G_{11} = 2m \,\psi^1_1\,\psi^1_1 +2m \,\psi^2_1\,\psi^2_1 +A
\,\psi^3_1\,\psi^3_1 +A \,\psi^4_1\,\psi^4_1
+(ma^2+2B) \,\psi^5_1\,\psi^5_1 \\
\phantom{G_{11}=}{} +ma\,\cos\psi\,\psi^1_1\,\psi^5_1
+ma\,\sin\psi\,\psi^2_1\,\psi^5_1
\\
\phantom{G_{11}}{}= 2m \,R^2\,\cos^2\psi +2m \,R^2\,\sin^2\psi +A
\,\frac{R^2}{a^2} +A  =2mR^2+A\left(1+\frac{R^2}{a^2} \right),
\\
G_{22} = 2m \,\psi^1_2\,\psi^1_2 +2m \,\psi^2_2\,\psi^2_2 +A
\,\psi^3_2\,\psi^3_2 +A \,\psi^4_2\,\psi^4_2
+(ma^2+2B) \,\psi^5_2\,\psi^5_2\\
\phantom{G_{22}=}{} +ma\,\cos\psi\,\psi^1_2\,\psi^5_2
+ma\,\sin\psi\,\psi^2_2\,\psi^5_2 = A \,\frac{R^2}{a^2}+ma^2+2B,
\\
G_{12} = 2m \,\psi^1_1\,\psi^1_2 +2m \,\psi^2_1\,\psi^2_2 +A
\,\psi^3_1\,\psi^3_2 +A \,\psi^4_1\,\psi^4_2
+(ma^2+2B) \,\psi^5_1\,\psi^5_2\\
\phantom{G_{12}=}{} +ma\,\cos\psi\,\psi^1_1\,\psi^5_2
+ma\,\sin\psi\,\psi^2_1\,\psi^5_2
\\
\phantom{G_{12}}{}= {}-A \frac{R^2}{a^2}
+ma\,\cos\psi\,R\,\cos\psi +ma\,\sin\psi\,R\,\sin\psi=maR-A
\frac{R^2}{a^2},
\end{gather*}
and we obtain
\[
[G_{\alpha\beta}]= \left[
\begin{array}{cc}
2mR^2+A\left(1+\frac{R^2}{a^2} \right) & maR-A \frac{R^2}{a^2}  \\[8pt]
maR-A \frac{R^2}{a^2} & A \,\frac{R^2}{a^2}+ma^2+2B
\end{array}
\right].
\]
It follows that
\begin{gather*}
G\doteq \det[G_{\alpha\beta}]
=\left[2mR^2+A\left(1+\frac{R^2}{a^2} \right)\right] \left[A
\,\frac{R^2}{a^2}+ma^2+2B \right]-\left[maR-A
\frac{R^2}{a^2}\right]^2
\\
\phantom{G}{}= 2mA \,\frac{R^4}{a^2}+2m^2R^2a^2+4mR^2B +A^2
\,\frac{R^2}{a^2}+ma^2A+2AB
+A^2\frac{R^4}{a^4}+ma^2A\frac{R^2}{a^2}\\
\phantom{G\doteq}{}+2AB\frac{R^2}{a^2} - m^2a^2R^2 -A^2
\frac{R^4}{a^4}+2mA \frac{R^3}{a}
\\
\phantom{G}{}=\frac{R^2}{a^2}\,
\left( 2mAR^2  + A^2 + ma^2A +2AB+2mAaR -m^2a^4 + 2m^2a^4 + 4ma^2B \right) \\
\phantom{G\doteq}{} +A\,(ma^2+2B)
\\
\phantom{G}{}=\frac{R^2}{a^2}\, \left( 2mAR^2  + A^2 + ma^2A
+2AB+2mAaR + m^2a^4  + 4ma^2B \right) +A\,(ma^2+2B)
\\
\phantom{G}{}=\frac{R^2}{a^2}\, \left( m\,(2AR^2+ a^2A +2AaR+
4a^2B+ ma^4)+ A^2 +2AB    \right) +A\,(ma^2+2B).
\end{gather*}
We observe that the determinant $G$ is a constant. The inverse
matrix is
\[
[G^{\alpha\beta}]=G^{-1} \left[
\begin{array}{cc}
A \,\frac{R^2}{a^2}+ma^2+2B & A \frac{R^2}{a^2} -maR \\[8pt]
A \frac{R^2}{a^2} -maR& 2mR^2+A\left(1+\frac{R^2}{a^2} \right)
\end{array}
\right].
\]
From the expression of the kinetic energy,
\[
K=m \,(\dot x^2+\dot y^2)+ \tfrac A2\,(\dot \theta_1^2+\dot
\theta_2^2)+ (\tfrac
12\,m\,a^2+B)\,\dot\psi^2+m\,a\,\dot\psi\;(\dot x\;\cos\psi+\dot
y\;\sin\psi),
\]
we obtain
 \[
\left[\dfrac{\partial K}{\partial \dot q^i}\right]= \left[
\begin{array}{c}
2m\dot x+ma\cos\psi \,\dot\psi \\
2m\dot y +ma\sin\psi\,\dot\psi \\
A\dot\theta_1 \\
A\dot\theta_2 \\
(ma^2+2B)\dot\psi+ma (\dot x\,\cos\psi+\dot y\,\sin\psi) \\
\end{array}
\right]
\]
and
\[
\left[\dfrac{\partial K}{\partial q^i}\right]= \left[
\begin{array}{c}
0 \\
0 \\
0 \\
0 \\
ma\dot\psi\,(\dot y\,\cos\psi-\dot x\,\sin\psi)\\
\end{array}
\right].
\]
The Lagrange equations for the free motions are
\begin{gather*}
2m\ddot x+ma\cos\psi \,\ddot\psi-ma\sin\psi \,\dot\psi^2= A_1,       \\
2m\ddot y +ma\sin\psi\,\ddot\psi +ma\cos\psi\,\dot\psi^2= A_2,    \\
A\ddot\theta_1=  A_3,     \\
A\ddot\theta_2=  A_4,    \\
(ma^2+2B)\ddot\psi+ma (\ddot x\cos\psi+\ddot y\sin\psi)-ma (\dot x\sin\psi-\dot y\,\cos\psi)\dot\psi \\
\qquad{}= ma\dot\psi(\dot y\cos\psi-\dot x\sin\psi) + A_5.
\end{gather*}
They show that
\begin{gather*}
L_1=A_1 +ma\sin\psi \,\dot\psi^2,      \\
L_2= A_2  -ma\cos\psi \,\dot\psi^2,  \\
L_3=A_3,     \\
L_4=A_4,    \\
L_5= A_5 +ma\dot\psi\,(\dot y\,\cos\psi-\dot x\,\sin\psi) +ma
(\dot x\,\sin\psi-\dot y\,\cos\psi)\,\dot\psi=A_5.
\end{gather*}
Thus, due to the parametric equations (\ref{2DPE}),
\begin{gather*}
\bar Z_1=A_1 +\tfrac {mR^2}a \, (z^2-z^1)^2\,\sin\psi,  \qquad
\bar Z_2= A_2-\tfrac {mR^2}a \, (z^2-z^1)^2\,\cos\psi,  \\
\bar Z_3=A_3, \qquad \bar Z_4=A_4, \qquad \bar Z_5= A_5.
\end{gather*}
Let us compute $Z_\alpha=\bar Z_i\,\psi^i_\alpha$ -- recall
(\ref{2DPSI}):
\begin{gather*}
Z_1= (A_1 +\tfrac {mR^2}a \, (z^2-z^1)^2\,\sin\psi)R\cos\psi +
(A_2-\tfrac {mR^2}a \, (z^2-z^1)^2\,\cos\psi) R\sin\psi \\
\phantom{Z_1=}{}-\tfrac Ra\,A_3+A_4
\\
\phantom{Z_1}{}=A_1\,R\,\cos\psi+A_2\,R\,\sin\psi-\tfrac
Ra\,A_3+A_4 =R\,(A_1\,\cos\psi+A_2\,\sin\psi-\tfrac 1a\,A_3)+A_4,
\\
Z_2=\tfrac Ra\,A_3+A_5.
\end{gather*}
Thus, the dynamical equations (\ref{ZZR}) associated with the
vector f\/ield $\mathbf Z$ are
\[
\begin{split}
&\dfrac{dx}{dt} = R\;\cos\psi\;z^1, \\
& \dfrac{dx}{dt}= R\;\sin\psi\;z^1, \\
&\dfrac{d\psi}{dt}=\frac Ra (z^2-z^1), \\
&\dfrac{d\theta_1}{dt} =z^1,\\
&\dfrac{d\theta_2}{dt}=z^2,
\end{split}
\qquad\quad
\begin{split}
& \dfrac{dz^1}{dt}=R\,(A_1\,\cos\psi+A_2\,\sin\psi-\tfrac 1a\,A_3)+A_4, \\
& \dfrac{dz^2}{dt}=\tfrac Ra\,A_3+A_5.
\end{split}
\]

\subsection{Two points with parallel velocities}

Two material points $P_1=(x_1,y_1)$ and $P_2=(x_2,y_2)$ running on
the Cartesian plane $\mathbb R^2=(x,y)$ are constrained to have
parallel  vector-velocities $\mathbf v_1$ and $\mathbf v_2$. This
is an example of non-linear non-holonomic constraint, since it is
expressed by the quadratic homogeneous equation
 \begin{equation}\label{PV}
C=\dot x_1\,\dot y_2-\dot x_2\,\dot y_1=0.
  \end{equation}
The conf\/iguration manifold is $Q_4=\mathbf R^4$ with ordered
Lagrangian coordinates
  \[
(q^1,q^2,q^3,q^4)=(x_1,y_1,x_2,y_2).
  \]
The kinetic energy is $K=\tfrac 12\, m_1\,(\dot x_1^2+\dot
y_1^2)+\tfrac 12\, m_2\,(\dot x_2^2+\dot y_2^2)$. Hence,
  \[
[g_{ij}]= \left[
\begin{array}{cccc}
m_1 & 0 & 0 & 0 \\
0 & m_1 & 0 & 0 \\
0 & 0 & m_2 & 0 \\
0 & 0 & 0 & m_2 \\
\end{array}
\right].
  \]
 (i) First method.
Since $\dim (Q)=4$ and $\dim (C)=7$, for a parametric
representation of the constraint $\dot x_1\,\dot y_2-\dot
x_2\,\dot y_1=0$  we need three parameters $(z^1,z^2,z^3)$. Let us
consider the parameters
$(z^\alpha)=(z^1,z^2,z^3)=(\rho,\sigma,\theta)$ and the parametric
equations
 \begin{equation}
  \begin{split}
& \dot x_1 = \rho \,\cos\theta ,\\
& \dot y_1 = \rho \,\sin \theta  ,
\end{split}
\qquad\quad
\begin{split}
& \dot x_2 = \sigma \,\cos\theta  ,\\
& \dot y_2 = \sigma \sin \theta  .
\end{split}
  \end{equation}
The meaning of the parameters is the following: $\rho^2=\mathbf
v_1^2=\dot x_1^2+\dot y_1^2$, $\sigma^2=\mathbf v_2^2=\dot
x_2^2+\dot y_2^2$, and $\theta$ is the  angle of the two vector
velocities w.r.to the $x$-axis. Then we f\/ind:
\begin{gather*}
[\psi^i_\alpha]= \left[
\begin{array}{cccc}
\cos\theta & \sin\theta & 0 & 0 \\
0 & 0 & \cos\theta & \sin\theta\\
{}-\rho\,\sin\theta & \rho\,\cos\theta & {}-\sigma\,\sin\theta & \sigma\,\cos\theta \\
\end{array}
\right] \qquad (\hbox{$\alpha$ index of line}),
\\
[G_{\alpha\beta}]\doteq[g_{ij}\,\psi^i_\alpha\psi^j_\beta]= \left[
\begin{array}{ccc}
m_1 & 0 & 0  \\
0 & m_2 & 0  \\
0 & 0 & m_1\rho^2+m_2\sigma^2  \\
\end{array}
\right],
\\
[G^{\alpha\beta}]= \left[
\begin{array}{ccc}
\dfrac 1{m_1} & 0 & 0  \\
0 & \dfrac 1{m_2} & 0  \\
0 & 0 & \dfrac 1{m_1\rho^2+m_2\sigma^2}  \\
\end{array}
\right].
\end{gather*}
Since,
  \[
\left[\dfrac{\partial K}{\partial \dot q^i}\right]=[m_1\dot x_1
\;,\; m_1\dot y_1   \;,\; m_2\dot x_2     \;,\; m_2\dot y_2
], \qquad \left[\dfrac{\partial K}{\partial q^i}\right]=[0\;,\;0
\;,\; 0 \;,\; 0],
 \]
the Lagrange equations for the free motions $g_{ij}\,\ddot
q^j=A_i-\Gamma_{hki}\,\dot q^h\,\dot q^k$ read
  \begin{equation}\label{2P1}
 m_1\,\ddot x_1=A_1  ,   \qquad
m_1\,\ddot y_1=A_2, \qquad m_2\,\ddot x_2=A_3 , \qquad m_2\,\ddot
y_2=A_4 .
  \end{equation}
They show that $\bar Z_i=A_i$. Hence,
\begin{gather*}
[Z_\alpha]\doteq [\psi^i_\alpha \bar Z_i]    =
\left[\begin{array}{c}
A_1\,\cos\theta+A_2\,\sin\theta \\
A_3\,\cos\theta+A_4\,\sin\theta \\
\rho\,(A_2\,\cos\theta-A_1\,\sin\theta)+\sigma\,(A_4\,\cos\theta-A_3\,\sin\theta)
\end{array}
\right],
\\
[Z^\alpha]\doteq[G^{\alpha\beta}\,Z_\beta]= \left[\begin{array}{c}
\dfrac {A_1\,\cos\theta+A_2\,\sin\theta}{m_1} \\
\dfrac{A_3\,\cos\theta+A_4\,\sin\theta}{m_2} \\
\dfrac{\rho\,(A_2\,\cos\theta-A_1\,\sin\theta)+\sigma\,(A_4\,\cos\theta-A_3\,\sin\theta)}{m_1\rho^2+m_2\sigma^2} \\
\end{array}
\right],
\end{gather*}
and the dif\/ferential system associated with $\mathbf Z$ is
\begin{equation}\label{Z2P}
  \begin{aligned}
& \dfrac {dx_1}{dt} = \rho \,\cos\theta ,\\
& \dfrac {dy_1}{dt} = \rho \,\sin \theta  ,\\
& \dfrac {dx_2}{dt} = \sigma \,\cos\theta  ,\\
& \dfrac {dy_2}{dt} = \sigma \sin \theta  ,
\end{aligned}
\qquad\qquad
  \begin{aligned}
& \dfrac{d\rho}{dt}=\dfrac {A_1\,\cos\theta+A_2\,\sin\theta}{m_1}, \\
& \dfrac{d\sigma}{dt}=\dfrac{A_3\,\cos\theta+A_4\,\sin\theta}{m_2}, \\
&
\dfrac{d\theta}{dt}=\dfrac{\rho\,(A_2\,\cos\theta-A_1\,\sin\theta)+
\sigma\,(A_4\,\cos\theta-A_3\,\sin\theta)}{m_1\rho^2+m_2\sigma^2},
\end{aligned}
  \end{equation}
where the Lagrangian active forces are in general known functions
of $(x_1,y_1,x_2,y_2)$ and $(\rho,\sigma,\theta)$.

In the special case of an inclined plane we have $A_1=m_1\,g$,
$A_3=m_2\,g$, $A_2=A_4=0$, and the system (\ref{Z2P}) becomes
\begin{equation}\label{Z2P1}
  \begin{aligned}
& \dfrac {dx_1}{dt} = \rho \,\cos\theta ,\\
& \dfrac {dy_1}{dt} = \rho \,\sin \theta  ,\\
& \dfrac {dx_2}{dt} = \sigma \,\cos\theta  ,\\
& \dfrac {dy_2}{dt} = \sigma \sin \theta  ,
\end{aligned}
\qquad\qquad
  \begin{aligned}
& \dfrac{d\rho}{dt}=g\,\cos\theta, \\
& \dfrac{d\sigma}{dt}=g\,\cos\theta, \\
&
\dfrac{d\theta}{dt}={}-g\,\sin\theta\,\dfrac{m_1\,\rho+m_2\,\sigma}
{m_1\rho^2+m_2\sigma^2}.
\end{aligned}
  \end{equation}
Note that the last three equations are separated from the f\/irst
four. This occurs in general when the Lagrangian active forces do
not depend on the position of the point, but only on their
velocities. For equal masses $m_1=m_2$, we have a further
simplif\/ication:
\begin{equation}\label{Z2P2}
  \begin{aligned}
& \dfrac {dx_1}{dt} = \rho \,\cos\theta ,\\
& \dfrac {dy_1}{dt} = \rho \,\sin \theta  ,\\
& \dfrac {dx_2}{dt} = \sigma \,\cos\theta  ,\\
& \dfrac {dy_2}{dt} = \sigma \sin \theta  ,
\end{aligned}
\qquad\qquad
  \begin{aligned}
& \dfrac{d\rho}{dt}=g\,\cos\theta, \\
& \dfrac{d\sigma}{dt}=g\,\cos\theta, \\
& \dfrac{d\theta}{dt}={}-g\,\sin\theta\,\dfrac{\rho+\sigma}
{\rho^2+\sigma^2}.
\end{aligned}
  \end{equation}

(ii) Second method for the single constraint equation (\ref{PV}).
In this case,
  \[
[C_i]= \left[\dot y_2\;,\;{}-\dot x_2 \;,\; {}-\dot y_1\;,\;\dot
x_1  \right]
  \]
does not have the maximal rank for $\mathbf v_1=\mathbf v_2=0$.
This is a singular state for whatever conf\/iguration (see Remark
\ref{r:p}); the set of the singular states is the
\emph{zero-section} of $TQ$. Moreover, since
\[
[g^{ij}]= \left[
\begin{array}{cccc}
\frac 1{m_1} & 0 & 0 & 0 \\
0 & \frac 1{m_1} & 0 & 0 \\
0 & 0 & \frac 1{m_2} & 0 \\
0 & 0 & 0 & \frac 1{m_2} \\
\end{array}
\right],
\]
we have:
\begin{gather*}
[C^i]\doteq [g^{ij}\,C_j]=\left[\frac{\dot
y_2}{m_1}\;,\;{}-\frac{\dot x_2}{m_1} \;,\; {}-\frac{\dot
y_1}{m_2}\;,\;\frac{\dot x_1}{m_2}  \right],
\\
G=\dfrac{\dot x_2^2+\dot y_2^2}{m_1}+\dfrac{\dot x_1^2+\dot
y_1^2}{m_2}=\dfrac{2K}{m_1\,m_2}, \qquad
G^{-1}=\dfrac{m_1\,m_2}{2K},
\\
[C_\ast^i]= G^{-1}\, \left[\frac{\dot y_2}{m_1}\;,\;{}-\frac{\dot
x_2}{m_1} \;,\; {}-\frac{\dot y_1}{m_2}\;,\;\frac{\dot x_1}{m_2}
\right],
\\
[\pi^{ij}]=G^{-1}\, \left[
\begin{array}{cccc}
\dfrac{\dot y_2^2}{m_1^2} & -\dfrac{\dot x_2\dot y_2}{m_1^2} &
-\dfrac{\dot y_1\dot y_2}{m_1\,m_2} & \dfrac{\dot x_1\dot
y_2}{m_1\,m_2}
\\[10pt]
-\dfrac{\dot x_2\dot y_2}{m_1^2} & \dfrac{\dot x_2^2}{m_1^2} &
\dfrac{\dot x_2\dot y_1}{m_1\,m_2} & -\dfrac{\dot x_1\dot
x_2}{m_1\,m_2}
 \\[10pt]
-\dfrac{\dot y_1\dot y_2}{m_1\,m_2} & \dfrac{\dot y_1\dot
x_2}{m_1\,m_2} & \dfrac{\dot y_1^2}{m_2^2} & -\dfrac{\dot x_1\dot
y_1}{m_2^2}
 \\[10pt]
\dfrac{\dot x_1\dot y_2}{m_1\,m_2} & -\dfrac{\dot x_1\dot
x_2}{m_1\,m_2} & -\dfrac{\dot x_1\dot y_1}{m_2^2}& \dfrac{\dot
x_1^2}{m_2^2}
\end{array}
\right],
\\
[g^{ij}-\pi^{ij}]= \left[
\begin{array}{cccc}
\dfrac 1{m_1} & 0 & 0 & 0 \\
0 & \dfrac 1{m_1} & 0 & 0 \\
0 & 0 & \dfrac 1{m_2} & 0 \\
0 & 0 & 0 & \dfrac 1{m_2} \\
\end{array}
\right]\\
\phantom{[g^{ij}-\pi^{ij}]=}{}- \dfrac{m_1\,m_2}{2K}\, \left[
\begin{array}{cccc}
\dfrac{\dot y_2^2}{m_1^2} & -\dfrac{\dot x_2\dot y_2}{m_1^2} & -\dfrac{\dot y_1\dot y_2}{m_1\,m_2} & \dfrac{\dot x_1\dot y_2}{m_1\,m_2} \\[10pt]
-\dfrac{\dot x_2\dot y_2}{m_1^2} & \dfrac{\dot x_2^2}{m_1^2} & \dfrac{\dot x_2\dot y_1}{m_1\,m_2} & -\dfrac{\dot x_1\dot x_2}{m_1\,m_2} \\[10pt]
-\dfrac{\dot y_1\dot y_2}{m_1\,m_2} & \dfrac{\dot y_1\dot x_2}{m_1\,m_2} & \dfrac{\dot y_1^2}{m_2^2} & -\dfrac{\dot x_1\dot y_1}{m_2^2} \\[10pt]
\dfrac{\dot x_1\dot y_2}{m_1\,m_2} & -\dfrac{\dot x_1\dot x_2}{m_1\,m_2} & -\dfrac{\dot x_1\dot y_1}{m_2^2}& \dfrac{\dot x_1^2}{m_2^2} \\
\end{array}
\right],
\\
g^{11}-\pi^{11}= \dfrac 1{m_1}-\dfrac{m_2}{2K}\dfrac{\dot
y_2^2}{m_1} =\dfrac 1{m_1}\left(1-\dfrac{m_2\,\dot
y_2^2}{2K}\right) =\dfrac{2K-m_2\,\dot y_2^2} {2\,m_1\,K}
=\dfrac{m_1\,(\dot x_1^2+\dot y_1^2)+ m_2\,\dot x_2^2}
{2\,m_1\,K},
\\
g^{12}-\pi^{12}= \dfrac{m_1\,m_2}{2K}\,\dfrac{\dot x_2\dot
y_2}{m_1^2} = \dfrac{m_2\,\dot x_2\dot y_2}{2\,m_1\,K},
\\
g^{13}-\pi^{13}= \dfrac{m_1\,m_2}{2K}\,\dfrac{\dot y_1\dot
y_2}{m_1\,m_2} =\dfrac{\dot y_1\dot y_2}{2K},
\\
g^{14}-\pi^{14}= {}-\dfrac{m_1\,m_2}{2K}\,\dfrac{\dot x_1\dot
y_2}{m_1\,m_2} ={}-\dfrac{\dot x_1\dot y_2}{2K},
\\
g^{22}-\pi^{22}= \dfrac 1{m_1}-\dfrac{m_2}{2K}\dfrac{\dot
x_2^2}{m_1} =\dfrac 1{m_1}\left(1-\dfrac{m_2\,\dot
x_2^2}{2K}\right) =\dfrac{2K-m_2\,\dot x_2^2} {2\,m_1\,K}
=\dfrac{m_1\,(\dot x_1^2+\dot y_1^2)+ m_2\,\dot y_2^2}
{2\,m_1\,K},
\\
g^{23}-\pi^{23}={} -\dfrac{m_1\,m_2}{2K}\,\dfrac{\dot x_2\dot
y_1}{m_1\,m_2} = {}-\dfrac{\dot x_2\dot y_1}{2K},
\\
g^{24}-\pi^{24}= \dfrac{m_1\,m_2}{2K}\,\dfrac{\dot x_1\dot
x_2}{m_1\,m_2} = \dfrac{\dot x_1\dot x_2}{2K},
\\
g^{33}-\pi^{33}=  \dfrac 1{m_2}-\dfrac{m_1}{2K}\dfrac{\dot
y_1^2}{m_2} =\dfrac 1{m_2}\left(1-\dfrac{m_1\,\dot
y_1^2}{2K}\right) =\dfrac{2K-m_1\,\dot y_1^2} {2\,m_2\,K}
=\dfrac{m_2\,(\dot x_2^2+\dot y_2^2)+ m_1\,\dot x_1^2}
{2\,m_2\,K},
\\
g^{34}-\pi^{34}=\dfrac{m_1\,m_2}{2K}\,\dfrac{\dot x_1\dot
y_1}{m_2^2} =\dfrac{m_1\,\dot x_1\dot y_1}{2\,m_2\,K} ,
\\
g^{44}-\pi^{44}=  \dfrac 1{m_2}-\dfrac{m_1}{2K}\dfrac{\dot
x_1^2}{m_2} =\dfrac 1{m_2}\left(1-\dfrac{m_1\,\dot
x_1^2}{2K}\right) =\dfrac{2K-m_1\,\dot x_1^2} {2\,m_2\,K}
=\dfrac{m_2\,(\dot x_2^2+\dot y_2^2)+ m_1\,\dot y_1^2}
{2\,m_2\,K},
\\
[\partial_iC]=[0 \;,\; 0 \;,\; 0 \;,\; 0], \qquad
[\partial_jC\,C_\ast^i]=[\mathbf 0].
\end{gather*}
As in this case $L_i=A_i$, we have
\begin{gather*}
D^i\doteq (g^{ij}-\pi^{ij})\,L_j-\dot
q^j\,\partial_jC\,C_\ast^i=(g^{ij}-\pi^{ij})\,L_j=(g^{ij}-\pi^{ij})\,A_j,
\\
D^1=\dfrac{1}{2K}\left(A_1\dfrac{m_1\,(\dot x_1^2+\dot y_1^2)+
m_2\,\dot x_2^2} {m_1} +A_2 \dfrac{m_2\,\dot x_2\dot y_2}{m_1}
+A_3\,\dot y_1\dot y_2 -A_4 \,\dot x_1\dot y_2\right),
\\
D^2=\dfrac{1}{2K}\left(A_1 \dfrac{m_2\,\dot x_2\dot y_2}{m_1}
+A_2\dfrac{m_1\,(\dot x_1^2+\dot y_1^2)+ m_2\,\dot y_2^2} {m_1}
-A_3\,\dot x_2\dot y_1 +A_4 \,\dot x_1\dot x_2\right),
\\
D^3=\dfrac{1}{2K}\left(A_1\dot y_1\dot y_2-A_2 \dot x_2\dot y_1
+A_3\dfrac{m_2\,(\dot x_2^2+\dot y_2^2)+ m_1\,\dot x_1^2} {m_2}
+A_4\dfrac{m_1\,\dot x_1\dot y_1}{m_2}\right),
\\
D^4= \dfrac{1}{2K}\left({}-A_1\dot x_1\dot y_2+A_2\dot x_1\dot x_2
+A_3\dfrac{m_1\,\dot x_1\dot y_1}{m_2} +A_4\dfrac{m_2\,(\dot
x_2^2+\dot y_2^2)+ m_1\,\dot y_1^2}{m_2}\right).
\end{gather*}
Then the dynamical system (\ref{DD}) reads
\begin{gather*}
\begin{split}
\dfrac{dx_1}{dt}=\dot x_1,
\\
\dfrac{dy_1}{dt}=\dot y_1,
\\
\dfrac{dx_2}{dt}= \dot x_2,
\\
\dfrac{dy_2}{dt}=\dot y_2,
\end{split}
\qquad
\begin{split}
&\dfrac{d\dot x_1}{dt}=\dfrac{1}{2K}\left(A_1\dfrac{m_1\,(\dot
x_1^2+\dot y_1^2)+ m_2\,\dot x_2^2} {m_1} +A_2 \dfrac{m_2\,\dot
x_2\dot y_2}{m_1} +A_3\,\dot y_1\dot y_2 -A_4 \,\dot x_1\dot
y_2\right),
\\
&\dfrac{d\dot y_1}{dt}=\dfrac{1}{2K}\left(A_1 \dfrac{m_2\,\dot
x_2\dot y_2}{m_1}     +A_2\dfrac{m_1\,(\dot x_1^2+\dot y_1^2)+
m_2\,\dot y_2^2} {m_1} -A_3\,\dot x_2\dot y_1 +A_4 \,\dot x_1\dot
x_2\right),
\\
&\dfrac{d\dot x_2}{dt}=\dfrac{1}{2K}\left(A_1\dot y_1\dot y_2-A_2
\dot x_2\dot y_1 +A_3\dfrac{m_2\,(\dot x_2^2+\dot y_2^2)+
m_1\,\dot x_1^2} {m_2} +A_4\dfrac{m_1\,\dot x_1\dot
y_1}{m_2}\right),
\\
&\dfrac{d\dot y_2}{dt}= \dfrac{1}{2K}\left({}-A_1\dot x_1\dot
y_2+A_2\dot x_1\dot x_2 +A_3\dfrac{m_1\,\dot x_1\dot y_1}{m_2}
+A_4\dfrac{m_2\,(\dot x_2^2+\dot y_2^2)+ m_1\,\dot
y_1^2}{m_2}\right).
\end{split}\!\!\!\!
\end{gather*}
For two points running on an inclined plane,
\[
\begin{split}
&D^1=g\,\dfrac {m_1\,(\dot x_1^2+\dot y_1^2)+ m_2\,(\dot x_2^2
+\dot y_1\dot y_2)}{m_1\,(\dot x_1^2+\dot y_1^2)+m_2\,(\dot
x_2^2+\dot y_2^2)},
\\
&D^2=g\,\dfrac {m_2\dot x_2\,(\dot y_2 -\dot y_1)}{m_1\,(\dot
x_1^2+\dot y_1^2)+m_2\,(\dot x_2^2+\dot y_2^2)},
\end{split}
\qquad
\begin{split}
&D^3=g\,\dfrac{m_2\,(\dot x_2^2+\dot y_2^2)+m_1\,(\dot x_1^2+\dot
y_1\dot y_2)}{ m_1\,(\dot x_1^2+\dot y_1^2)+m_2\,(\dot x_2^2+\dot
y_2^2)},
\\
&D^4= g\, \dfrac {m_1 \,\dot x_1\,(\dot y_1 -\dot y_2)}{
m_1\,(\dot x_1^2+\dot y_1^2)+m_2\,(\dot x_2^2+\dot y_2^2)}.
\end{split}
\]
For equal masses, $m_1=m_2$,
\[
\begin{split}
&D^1=g\,\dfrac {\dot x_1^2+\dot y_1^2+\dot x_2^2 +\dot y_1\dot
y_2}{\dot x_1^2+\dot y_1^2+\dot x_2^2+\dot y_2^2},
\\
&D^2=g\,\dfrac {\dot x_2\,(\dot y_2 -\dot y_1)}{\dot x_1^2+\dot
y_1^2+\dot x_2^2+\dot y_2^2},
\end{split}
\qquad
\begin{split}
&D^3=g\,\dfrac{\dot x_2^2+\dot y_2^2+\dot x_1^2+\dot y_1\dot
y_2}{\dot x_1^2+\dot y_1^2+\dot x_2^2+\dot y_2^2},
\\
&D^4= g\, \dfrac {\dot x_1\,(\dot y_1 -\dot y_2)}{\dot x_1^2+\dot
y_1^2+\dot x_2^2+\dot y_2^2}.
\end{split}
\]

\begin{note}
It is easy (and obvious) to propose examples of non-linear
constraints: it is suf\/f\/icient to choose any set of non-linear
independent equations $C^a(q,\dot q)=0$. However, any example of a
non-linear constraint remains meaningless unless we know how to
realize it \emph{physically} by means of realizable  devices. The
famous  Appell--Hamel example gives a matter of discussion (see
\cite{Neimark}, Ch.~4, \S~2). Indeed, in order to be really a
non-linear device, a certain distance of the Appell--Hamel device
must be inf\/initesimally small. This f\/its with the thought of
Hertz: non-linear constraints can be realized by passing to the
limit $x\to 0$ of certain physical quantities~$x$ (masses,
lengths, etc.) in devices realizing linear constraints.

The same kind of problem arises in trying to `\emph{realize}' two
mass-points moving with parallel velocities. A tentative project
has been presented in \cite{Benenti}. In fact, for an ef\/fective
project, we have to invent devices for
\begin{enumerate}\itemsep=0pt
\item Realizing a mass-point. \item Realizing a parallel transport
on the plane. \item Transforming the direction of the velocity of
a point into a solid segment. \item Applying forces of special
kind to the points (the weight is of course always present).
\end{enumerate}
This research is a work in progress. Updated information will be
found on my personal web-page.
\end{note}

\subsection*{Acknowledgments}
A preliminary version of this paper has been elaborated and
exposed at the University of Link\"oping, Department of
Mathematics, on May 27, 2005. I wish to thank Stefan Rauch and all
the Link\"oping school for their warm hospitality. I wish also to
thank: Waldyr Oliva, Willy Sarlet, David Martin de Diego for
making me aware of their contributions to the theory; Enrico
Pagani and Enrico Massa, for the enlightening discussions during a
workshop on dynamical systems held in Torino in April 2005 --
their papers \cite{Massa-Pagani-1991, Massa-Pagani-1997} have been
of great help; Beppe Gaeta, for pointing me out some errors in my
f\/irst manuscript; Franco Cardin for his kind invitation to give
a seminar in Padova on the contents of this paper (November 2006).

\pdfbookmark[1]{References}{ref}
\LastPageEnding

\end{document}